\documentclass[a4paper, reqno]{amsart}
\usepackage{amssymb, a4wide, amsmath, amsfonts}
\usepackage[usenames, dvipsnames]{xcolor}
\usepackage{enumitem}
\usepackage{mathtools}
\usepackage{thmtools}
\usepackage{bm}
\usepackage{mathrsfs}
\usepackage{xfrac}
\usepackage{bm}

\theoremstyle{plain}
\declaretheorem[name=Proposition, numberwithin=section]{proposition}
\declaretheorem[name=Theorem, sibling=proposition]{theorem}
\declaretheorem[name=Lemma, sibling=proposition]{lemma}

\declaretheorem[name=Definition, sibling=proposition]{definition}

\declaretheorem[name=Example, sibling=proposition]{example}

\usepackage[colorlinks, linkcolor=blue]{hyperref}

\newcommand\numberthis{\addtocounter{equation}{1}\tag{\theequation}}
\newcommand{\vertiii}[3]{{\vert\kern-0.25ex\vert\kern-0.25ex\vert #1 
		\vert\kern-0.25ex\vert\kern-0.25ex\vert_{#2}^{#3}}}
	
\newcommand{\norm}[2]{{\lVert #1\rVert_{#2}}}
\newcommand{\normto}[3]{{\lVert #1\rVert_{#2}^{#3}}}

\newcommand{\absh}[2]{{\left\lvert #1\right\rvert_{#2}}}

\newcommand{\omg}[2]{\omega_{\Phi_\beta}(#1, #2)}
\newcommand{\omgd}[3]{\omega^{#3}_{\Phi_\beta}(#1, #2)}

\newcommand{\lsmT}[1]{\lesssim_{#1}^T}
\newcommand{\lsmM}[1]{\lesssim_{#1}^M}
\newcommand{\lsmTM}[1]{\lesssim_{#1}^{T, M}}

\newcommand{\lsmToneMzero}[1]{\lesssim_{#1}^{T_1, M_0}}

\NewDocumentCommand{\mycommand}{ m O{} }{%
	Argument 1: #1%
	\IfValueT{#2}{, Argument 2: #2}%
	\IfNoValueT{#2}{, Argument 2: DefaultValue}%
}

\newcommand{\dd}{\,\mathrm{d}}

\newcommand{\Lcal}{\mathcal{L}}
\newcommand{\V}{\mathbb{V}}

\newcommand{\I}{\mathcal{I}}
\newcommand{\Iscr}{\mathscr{I}}
\newcommand{\X}{\mathbf{X}}
\newcommand{\R}{\mathbb{R}}
\newcommand{\N}{\mathbb{N}}
\newcommand{\Xbb}{\mathbb{X}}

\newcommand{\Xito}{\mathbf{X}^{\mathrm{It\hat{o}}}}
\newcommand{\Xbbito}{\mathbb{X}^{\mathrm{It\hat{o}}}}
\newcommand{\Xstrat}{\mathbf{X}^{\mathrm{Strat}}}
\newcommand{\Xbbstrat}{\mathbb{X}^{\mathrm{Strat}}}

\newcommand{\BOone}{B^\alpha_{\Phi_\beta, q}}
\newcommand{\BOtwo}{B^{2\alpha; 2}_{\Phi_{\beta/2}, q/2}}
\newcommand{\BOthree}{B^{3\alpha; 2}_{\Phi_{\beta/3}, q/3}}
\newcommand{\BOthreer}{B^{3\alpha; 2}_{\Phi_{\beta/3}, \bar{q}}}
\newcommand{\BOgamma}{B^{\gamma; 2}_{\Phi_{\beta/2}, q/2}}
\newcommand{\BOd}[1]{B^{\alpha; #1}_{\Phi_{\beta}, q}}
\newcommand{\BOrp}{\mathscr{B}^\alpha_{\Phi_\beta, q}}
\newcommand{\BO}[3]{B^{#1}_{\Phi_{#2}, #3}}

\newcommand{\BOSone}{B^{\alpha; 2}_{\Phi_{\beta_1}, q_1}}
\newcommand{\BOStwo}{B^{\gamma; 2}_{\Phi_{\beta_2}, q_2}}
\newcommand{\BOSthree}{B^{\gamma; 3}_{\Phi_{\beta_2}, q_2}}
\newcommand{\BOSr}{B^{\gamma; 2}_{\Phi_{\beta_2}, r}}

\newcommand{\OSletter}{L}
\newcommand{\Osp}{\OSletter^{\Phi_\beta}}
\newcommand{\Ospb}[1]{\OSletter^{\Phi_{#1}}}
\newcommand{\OspF}{\OSletter_{\mathrm{fin}}^{\Phi_\beta}}

\newcommand{\Lq}{L^q([0, T],\frac{\dd\tau}{\tau})}

\newcommand{\onT}{([0,T])}

\newcommand{\Rmk}{\R^{m\times k}}
\newcommand{\Rmn}{\R^{m\times n}}
\newcommand{\Rborpcodo}{\R^n}

\newcommand{\YscrT}{\mathscr{Y}_T}
\newcommand{\YscrTzero}{\mathscr{Y}_{T_0}}
\newcommand{\YscrTone}{\mathscr{Y}_{T_1}}
\newcommand{\Yn}{Y^{(n)}}
\newcommand{\Ynp}{Y'^{(n)}}
\newcommand{\BOepsd}[1]{\BO{#1}{\beta/2}{q/2}}
\newcommand{\BOeps}{\BOepsd{\alpha-\varepsilon;2}}
\newcommand{\Zscr}{\mathscr{Z}}

\newcommand{\CRPs}{\mathscr{D}^{2\alpha; X}_{\Phi_{\beta / 2}, q / 2}}
\newcommand{\CRPtld}{\mathscr{D}^{2\alpha; \tilde{X}}_{\Phi_{\beta / 2}, q / 2}}
\newcommand{\dX}{d_{X, \mathscr{D}^{2\alpha}_{\Phi_{\beta / 2}, q / 2}}}
\newcommand{\dXX}{d_{X, \tilde{X}, \mathscr{D}^{2\alpha}_{\Phi_{\beta / 2}, q / 2}}}
\newcommand{\dXon}[1]{d_{X, \mathscr{D}^{2\alpha}_{\Phi_{\beta / 2}, q / 2}([0,#1])}}
\newcommand{\dXXon}[1]{d_{X, \tilde{X}, \mathscr{D}^{2\alpha}_{\Phi_{\beta / 2}, q / 2}([0,#1])}}

\newcommand{\Bone}{B^\alpha_{p, q}}
\newcommand{\Besov}[3]{B^{#1}_{#2, #3}}

\newcommand{\BYMM}[1]{B_{\YscrT}((\bar{Y},\bar{Y}'),#1)}
\newcommand{\BYM}{\BYMM{M}}
\newcommand{\BYMMT}[2]{B_{#2}((\bar{Y},\bar{Y}'),#1)}

\newcommand{\hook}{\hookrightarrow}

\makeatletter
\@namedef{subjclassname@2020}{\textup{2020} Mathematics Subject Classification}
\makeatother

\begin{document}

\title[Differential equations driven by Besov-Orlicz paths]{Differential equations driven by Besov-Orlicz paths}

\author{Petr \v{C}oupek}
\address{Charles University\\
	Faculty of Mathematics and Physics\\
	Sokolovsk\'{a}~83\\
	186~75\\
	Prague~8\\
	Czech Republic}
\email[Corresponding author]{coupek@karlin.mff.cuni.cz}

\author{Franti\v{s}ek Hendrych}
\address{Charles University\\
	Faculty of Mathematics and Physics\\
	Sokolovsk\'{a}~83\\
	186~75\\
	Prague~8\\
	Czech Republic}
\email{hendrychfrantisek@karlin.mff.cuni.cz}

\author{Jakub Slav\'ik}
\address{Czech Academy of Sciences\\
	Institute of Information Theory and Automation\\
	Pod Vod\'{a}renskou v\v{e}\v{z}\'{i}~4\\
	180~00\\
	Prague~8\\
	Czech Republic}
\email{slavik@utia.cas.cz}



\begin{abstract}
	In the article, the rough path theory is extended to cover paths from the exponential Besov-Orlicz space \[B^\alpha_{\Phi_\beta,q}\quad\mbox{ for }\quad \alpha\in (1/3,1/2],\,\quad  \Phi_\beta(x) \sim \mathrm{e}^{x^\beta}-1\quad\mbox{with}\quad \beta\in (0,\infty), \quad\mbox{and}\quad q\in (0,\infty],\] and the extension is used to treat nonlinear differential equations driven by such paths. The exponential Besov-Orlicz-type spaces, rough paths, and controlled rough paths are defined and analyzed, a sewing lemma for such paths is given, and the existence and uniqueness of the solution to differential equations driven by these paths is proved. The results cover equations driven by paths of continuous local martingales with Lipschitz continuous quadratic variation (e.g.\ the Wiener process) or by paths of fractionally filtered Hermite processes in the $n$\textsuperscript{th} Wiener chaos with Hurst parameter $H\in (1/3,1/2]$ (e.g.\ the fractional Brownian motion).
\end{abstract}

\maketitle

\section{Introduction}
In the last two decades, rough path theory and its extensions have had enormous impact in the field of differential equations (DEs) driven by singular functions. Such singularity was originally described in terms of H\"older continuity (or, almost equivalently, in terms of finite $p$-variation) but in recent years, several extensions to paths of Sobolev or Besov regularity have been given (see the series of papers \cite{LiuProTei21,LiuProTei23,LiuProTei23a} and \cite{FriPro18,FriSee22, ProTra16}, respectively).

From the perspective of stochastic differential equations (SDEs), i.e.\ DEs driven by paths generated from a stochastic processes, such extensions are very useful. Indeed, if we consider these equations, one would expect that the solution retains the noise regularity because, roughly speaking, the solution should behave like the noise on small scales. Consider, for example, the Wiener process $W=(W_t, t\in [0,1])$ as the noise source. By a straightforward application of the Kolmogorov continuity theorem, it is immediately seen that its paths lie, almost surely, in the H\"older space $C^{\frac{1}{2}-\varepsilon}([0,1])$ for any $\varepsilon>0$. One can therefore fix such $\varepsilon>0$ and use the rough path machinery for H\"older continuous functions (e.g.\ \cite[Theorem 8.4]{FriHai20}) to enhance the Wiener path $W(\omega)$ to a Wiener rough path $\mathbf{W}(\omega)=(W(\omega),\mathbb{W}(\omega))$, where $\mathbb{W}(\omega)$ is the corresponding path of the, say, Stratonovich integral 
\[
    \mathbb{W}_{s,t}(\omega) = \int_s^t (W_r-W_s)\circ\dd{W}_r (\omega) = \frac{1}{2} (W_t(\omega)-W_s(\omega))^2,
\] 
to obtain the global solution to the rough differential equation (RDE)
	\begin{equation}
	\label{eq:SDE}
	\dd{Y}_t(\omega) = f(Y_t(\omega))\dd\mathbf{W}_t(\omega), \quad Y_0(\omega) = y,
	\end{equation}
for $y\in\R$ and $f\in C^3_b(\R)$ that will again be of $C^{1/2-\varepsilon}$-regularity. The solution obtained in this manner then agrees with the solution to the corresponding Stratonovich SDE. However, somehow one feels that more information about the solution can be obtained if we could employ some extra information about the regularity of the Wiener path. It is known, for example, that Wiener paths belong, almost surely, to the Besov space $B_{p,\infty}^{1/2}([0,1])$ for all $p\in [1,\infty)$ but not to the space $B_{p,q}^{1/2}([0,1])$ for any $q<\infty$ (see \cite{Cie91}). Therefore, by choosing $p\in [1,\infty)$ and by appealing to the Besov extension of rough paths in \cite{FriSee22} (namely to Theorem 5.6 therein), one in fact obtains that the solution to the RDE \eqref{eq:SDE} is actually of $B_{p,\infty}^{1/2}$-regularity. On the other hand, it is also known that Wiener paths belong, almost surely, to the modular H\"older space $C^{|r\log r|^{1/2}}([0,1])$, see \cite{Lev37}, but this space and the Besov space $B_{p,\infty}^{1/2}([0,1])$ are not included in one another. It was then soon realized that one can quantify the asymptotic growth of the $L^p$-modulus of continuity inside the $B_{p,\infty}^{1/2}$-norm and show that Wiener paths belong, almost surely, to the exponential Besov-Orlicz space $B_{\Phi_2,\infty}^{1/2}([0,1])$, where $\Phi_2(x)= \mathrm{e}^{x^2}-1$, see \cite{Cie93}, that lies in the intersection of the two spaces.

In fact, the exponential Besov-Orlicz spaces \[B_{\Phi_\beta,\infty}^\alpha([0,1])\quad\mbox{for}\quad \alpha\in (0,1)\quad \mbox{and} \quad \Phi_{\beta}\sim \mathrm{e}^{x^\beta}-1 \quad\mbox{with}\quad  \beta\in (0,\infty)\] appear to form a very natural scale of function spaces for a multitude of stochastic processes. For example, such path regularity is obtained for continuous local martingales with Lipschitz continuous quadratic variation (whose prototypical example is the Wiener process) or for the fractionally filtered Hermite processes in the $n$\textsuperscript{th} Wiener chaos with the Hurst parameter $H\in (0,1)$ \cite{BaiTaq14} (with examples such as the fractional Brownian motion \cite{DecUst99} or the Rosenblatt process \cite{Tud08}) and while the former processes are known to have paths in $B_{\Phi_2,\infty}^{1/2}([0,1])$, see \cite[Theorem~4.1]{OndSimKup18}, the latter have paths in $B_{\Phi_{2/n},\infty}^{H}([0,1])$, see \cite[Corollary 4.2]{CouOnd24}. Other results on Besov-Orlicz regularity of stochastic processes are also given in \cite{OndVer20,Wich21}.

In the present article, we therefore aim to solve differential equations of the form
	\[ \dd{Y}_t = f(Y_t)\dd{X}_t, \quad Y_0=y, \] 
on the interval $[0,T]$ for $T\in (0,\infty)$, $y\in\mathbb{R}^m$, $f\in C^3_b(\mathbb{R}^m, \mathbb{R}^{m\times n})$, and a path $X$ in the exponential Besov-Orlicz space $B_{\Phi_\beta,q}^\alpha([0,T];\mathbb{R}^n)$, $m,n\in\mathbb{N}$. In order to do so, we employ the rough path machinery; that is, initially, we extend the definition of the classical exponential Besov-Orlicz spaces to the exponential Besov-Orlicz-type spaces of multivariate maps that are suitable for rough path analysis and give several of their properties. This is done in \autoref{sec:preliminaries_BO}. We then proceed, in \autoref{sec:rps}, with the definition of an exponential Besov-Orlicz rough path $\X$ and show that the paths of the stochastic processes mentioned above can be indeed lifted to such rough paths. Subsequently, controlled rough paths are defined and their properties such as their stability under compositions with $C^2_b$-functions are given. We then prove a sewing lemma that is subsequently used to define a rough integral for paths of the considered Besov-Orlicz regularity. Finally, in \autoref{sec:RDEs}, we consider the (rough) DE
	\[ \dd{Y}_t = f(Y_t)\dd\X_t,\quad Y_0=y,\]
and we give the main result of the article in \autoref{thm:existanduniq} where we show that the equation admits a unique solution of $B_{\Phi_\beta,q}^\alpha$-regularity. As a consequence, not only do we obtain $B_{\Phi_2,\infty}^{1/2}$-regularity of the solution to equation \eqref{eq:SDE}, which improves the known results on DEs driven by Wiener paths, but we also obtain Besov-Orlicz regularity of Stratonovich-type DEs driven by paths of other, possibly non-Gaussian and non-Markovian, stochastic processes.

\section{Preliminaries: Hölder and Besov spaces}
\label{sec:preliminaries}

Let us begin by listing the basic notation and function spaces used throughout the article and by recalling the definitions of H\"older-type and Besov-type spaces used in the theory of rough paths.

\subsection{Basic notation and function spaces}
\label{sec:notation}

We use the following convention throughout the article: We write $A\lesssim B$ if there exists a finite positive constant $C$ such that $A\leq CB$. If $C$ depends on some  parameter $\theta$, we write either $A \lesssim_\theta B$ or $A \leq C(\theta) B$. If $A$, $B$, and $C$ also depend on an additional parameter $T>0$, i.e.\ the inequality $A(T)\leq C(\theta,T) B(T)$ holds for all $T>0$, we may wish to stress that $C(\theta, \cdot)$ is nondecreasing (in particular, the value of $C(\theta,T)$ will not tend to infinity as $T\to 0+$), we write $A(T)\lsmT{\theta}B(T)$. The value of the constant itself can change from one line to another without any additional comment.

In the list of function spaces below, let $n \in \N$, let $(\V_1,\lvert\,\cdot\,\rvert_{\V_1})$ and $(\V_2,\lvert\,\cdot\,\rvert_{\V_2})$ be finite-dimensional normed vector spaces, let $(E, d)$ be a nonempty metric space and let $V \subseteq \V_1$. We note that the definitions of the spaces of $E$-valued functions below formally depend on the particular choice of $e_0 \in E$. However, with a~different choice of $e_0$, only the respective norms differ while the function spaces themselves remain the same. For simplicity, we therefore assume that $e_0 \in E$ is fixed. When $E$ is a Banach space, we naturally choose $e_0=0$.

\begin{itemize}
	\item For $m, k \in \N$, we denote the space of $(m\times k)$-matrices by $\Rmk$. With $b \in \R^m$ and $\tilde{b} \in \R^k$, we note that the Hilbert–Schmidt norm $\lvert\,\cdot\,\rvert_{\R^{m\times k}}$ is compatible, i.e.\ $\lvert b\tilde{b}^\top\rvert_{\R^{m\times k}}\leq\lvert b\rvert_{\R^m}\lvert \tilde{b}\rvert_{\R^k}$, and symmetric, i.e.\ $\lvert b\tilde{b}^\top\rvert_{\R^{m\times k}}=\lvert \tilde{b}b^\top\rvert_{\R^{k\times m}}$.
	\item We denote the space of all bounded linear operators from $\V_1$ to $\V_2$ by $\mathcal{L}(\V_1; \V_2)$ and the space of $n$-ary linear operator from the $n$-product space $\V_1 \times \V_1 \times \dots \times \V_1$ to $\V_2$ by $\mathcal{L}^{(n)}(\V_1 \times \V_1 \times \dots \times \V_1; \V_2)$.
	\item We denote the space of all continuous functions from $V$ to $E$ by $C(V;E)$. Similarly, the space of all bounded continuous functions from $V$ to $E$ is denoted by $C_b(V;E)$. We equip $C_b(V; E)$ with the supremum norm $\norm{f}{C_b(V;E)} = \sup_{v \in V} d(f_v,e_0)$.
	\item If $V$ is open, $f \in C(V;\V_2)$, and $f$ is Fr\'echet-differentiable in $V$, we denote the Fr\'echet derivative of $f$ by $Df: V \to \mathcal{L}(\V_1;\V_2)$. Similarly, the symbol $D^n f: V \to \mathcal{L}^{(n)}(\V_1 \times \V_1 \times \dots \times \V_1;\V_2)$ stands for the $n$-th Fr\'echet derivative if it exists.
	\item For $n \in \N$, we denote

	\[
		\quad\quad\quad C_b^n(V; \V_2) = \left\{ f: V \to \V_2 \, \left| \, D^k f\in C_b(V;\mathcal{L}^{(k)}(\V_1 \times \V_1 \times \dots \times \V_1;\V_2)), k=0,\dots,n \right.\right\},
	\]
 
	\noindent where, for $k=0$, we identify the respective space of $0$-ary operators with $\V_2$. We equip the space $C_b^n(V;\V_2)$ with the norm 
 
    \begin{minipage}{\linewidth}
    \[
        \norm{f}{C_b^n(V;\V_2)}=\sum_{k=0}^n \norm{D^kf}{C_b(V; \mathcal{L}^{(k)}(\V_1 \times \V_1 \times \dots \times \V_1;\V_2))}.
    \]
    \end{minipage}
    
	\item For $\alpha \in (0, 1)$, the $\alpha$-H\"older continuous functions $C^\alpha(V;E)$ are defined by 

    \begin{minipage}{\linewidth}
    \[
        C^\alpha(V;E) = \{ f \in C(V;E) \mid [f]_{C^\alpha(V;E)} < \infty \},
    \]
    \end{minipage}
    
 \noindent where
 
    \begin{minipage}{\linewidth}
	\[
		[f]_{C^\alpha(V;E)} = \sup\limits_{\substack{v, \tilde{v} \in V\\ v \neq \tilde{v}}} \frac{d(f_v,f_{\tilde{v}})}{\lvert v-\tilde{v}\rvert_{\V_1}^\alpha}.
	\]
    \end{minipage}
    
	\item For $\mathcal{O} \subseteq \R^n$, $\mathfrak{L}(\mathcal{O})$ denotes the $\sigma$-algebra of Lebesgue measurable subsets of $\mathcal O$ and~$\mathfrak{B}(E)$ denotes the $\sigma$-algebra of Borel measurable subsets of $E$. We write $L^0(\mathcal O;E)$ for the set of equivalence classes of measurable functions $f:(\mathcal O, \mathfrak{L}(\mathcal O))\to (E, \mathfrak{B}(E))$ with respect to equality almost everywhere.
	\item For $p \in (0, \infty]$ and $f \in L^p([0, T]; E)$, we define the $L^p$-modulus of continuity by 
 
    \begin{minipage}{\linewidth}
    \[
        \omega_p(f, \tau) = \sup_{h \in [0, \tau]} \|d(f_\cdot,f_{\cdot+h})\|_{L^p([0, T-h])}, \quad \tau \in [0, T].
    \]
    \end{minipage}
    
	\item Let $\alpha \in (0, 1)$ and $p, q \in (0, \infty]$. The Besov space $B^\alpha_{p,q}([0, T];E)$ is defined by
	
    \begin{minipage}{\linewidth}
    \[
		B^\alpha_{p,q}([0, T]; E) = \left\{ f \in L^p([0, T]; E) \, \left| \, [f]_{B^\alpha_{p,q}([0, T]; E)} < \infty \right.\right\},
	\]
    \end{minipage}
 
	\noindent where

    \begin{minipage}{\linewidth}
	\[
		[f]_{B^\alpha_{p,q}([0, T];E)} = \left\|\frac{\omega_p(f, \tau)}{\tau^\alpha} \right\|_{L^q([0,T],\frac{\dd\tau}{\tau})}.
	\]
    \end{minipage}
\end{itemize}

\subsection{H\"older-type and Besov-type spaces} We now recall the definitions of H\"older-type and Besov-type spaces used in the theory of rough paths. The Besov-type spaces were first defined in \cite{FriSee22} to which we refer the reader for a more detailed exposition. Let us fix $d \in \{ 2, 3 \}$, $T \in (0, \infty)$, and a normed vector space $(\V, \lvert\,\cdot\,\rvert_{\V})$ for the rest of this section. We denote
\[
\triangle^d[0,T]=\{(u_1,\dots, u_d)\in[0,T]^d\,|\,u_1\leq\dots\leq u_d\}.
\]
For $f: [0, T] \to \V$ and $\varXi: \triangle^2[0, T] \to \V$, we define $\delta f: \triangle^2[0, T] \to \V$ and $\delta \varXi: \triangle^3[0, T] \to \V$ by
\begin{alignat*}{2}
 & \delta f_{s,t}  = f_t - f_s,   && (s,t)\in\triangle^2[0,T],\\ & \delta \varXi_{s,u,t} = \varXi_{s,t}-\varXi_{s,u}-\varXi_{u,t},  && \quad (s,u,t) \in \triangle^3[0, T].
\end{alignat*}
We note that $\delta(\delta f)\equiv 0$. Next, we construct the space of measurable functions on the simplex $\triangle^d[0, T]$. Let $\varXi, \tilde{\varXi}: (\triangle^d[0,T], \mathfrak{L}(\triangle^d[0,T]))\to (\V, \mathfrak{B}(\V))$. For $d=2$, we define
\[
	\varXi \sim_2 \tilde{\varXi}\ \text{if} \ \varXi_{r,r+h} = \tilde{\varXi}_{r,r+h} \ \text{for all $h \in [0, T]$ and almost all $r \in [0, T-h]$},
\]
and, similarly, for $d=3$, let
\[
	\varXi \sim_3 \tilde{\varXi}\ \text{if} \ \varXi_{r,r+\theta h,r+h} = \tilde{\varXi}_{r,r+\theta h,r+h} \ \text{for all $h \in [0, T]$, all $\theta \in [0, 1]$, and almost all $r \in [0, T-h]$}.
\]
It is straightforward to show that $\sim_d$ is an equivalence. We then define $L^{0;d}([0,T];\V)$ as the space of all equivalence classes of measurable functions $(\triangle^d[0,T], \mathfrak{L}(\triangle^d[0,T]))\to (\V, \mathfrak{B}(\V))$ with respect to equivalence $\sim_d$.

We can now recall the definition of the H\"older and Besov-type spaces. For $\alpha\in(0,\infty)$, we define the $\alpha$-H\"older-type space by
\[
	C^{\alpha; 2}([0,T];\V) = \left\{ \varXi: \triangle^2[0, T] \to \V \, \left| \, \| \varXi \|_{C^{\alpha; 2}([0,T]; \V)} < \infty \right.\right\},
\]
where
\[
	\| \varXi \|_{C^{\alpha; 2}([0,T];\V)} = \sup\limits_{\substack{0 < h \leq T\\ 0 \leq r \leq T-h}} \frac{ \absh{\varXi_{r, r+h}}{\V}}{|h|^\alpha}.
\]
Note that if $\alpha\in(0,1)$, then $[f]_{C^\alpha([0, T];\V)}=\norm{\delta f}{C^{\alpha;2}([0, T];\V)}$.

For $p \in (0, \infty]$ and $\varXi\in L^{0;d}([0,T];\V)$, the $L^p$-modulus of continuity of $\varXi$ is defined by
\begin{equation*}
	\omega^{d}_p(\varXi, \tau) = \begin{cases}
			\displaystyle \sup_{h \in [0, \tau]} \| \varXi_{\cdot, \cdot+h} \|_{L^p([0, T-h];\V)},& \quad d=2,\\
			\displaystyle \sup_{\theta\in[0, 1]}\sup_{h\in[0, \tau]}\norm{\varXi_{\cdot, \cdot+\theta h, \cdot+h}}{L^p([0, T-h];\V)},& \quad d=3,
		\end{cases}
        \quad\tau\in[0,T],
\end{equation*}
and for $\alpha\in(0,\infty)$ and $p,q \in (0, \infty]$, we define the Besov-type space $B^{\alpha; d}_{p,q}([0, T];\V)$ by
\[
	B^{\alpha; d}_{p,q}([0, T];\V) = \left\{ \varXi \in L^{0;d}([0,T];\V) \, \left| \, \norm{\varXi}{B^{\alpha; d}_{p,q}([0, T];\V)} < \infty \right.\right\},
\]
where
\[
	\norm{\varXi}{B^{\alpha; d}_{p,q}([0, T];\V)} = \left\lVert\frac{\omega_p^d(\varXi, \tau)}{\tau^\alpha} \right\rVert_{L^q([0,T];\frac{\dd\tau}{\tau})}.
\]
Note that for $\alpha \in (0, 1)$, it holds that $[f]_{B^\alpha_{p,q}([0, T]; \V)}=\norm{\delta f}{B^{\alpha; 2}_{p,q}([0, T]; \V)}$.

\subsection{Exponential Orlicz spaces}

Let us first recall some basic facts on exponential Orlicz spaces that will be needed for our analysis. For a thorough exposition, we refer the reader to the excellent monographs \cite{PicKufJohFuc13} and \cite{RaoRen91}. Let $\beta\in(0, \infty)$ and set $x_\beta=(\frac{1-\beta}{\beta})^{1/\beta}$ for $\beta<1$ and $x_\beta=0$ for $\beta\geq 1$. Define $\Psi_\beta, E_\beta, \Phi_\beta: [0,\infty) \to [0,\infty)$ by 
\begin{align*}
	\Psi_\beta(x) & = \exp({x^\beta}) - 1,
	\\
	E_\beta(x) & = 
	\begin{cases}
		\displaystyle \Psi_\beta(x), & x\in[x_\beta, \infty),\\
		\displaystyle \Psi_\beta(x_\beta) + \Psi_\beta'(x_\beta)(x - x_\beta), & x\in[0, x_\beta),
	\end{cases}
	\\
	\Phi_\beta(x) & = E_\beta(x) - E_\beta(0),
\end{align*}
for $x\geq 0$. It is easily seen that $\Phi_\beta$ is a convex function such that $\Phi_\beta(0) = 0$ and $\lim\limits_{x\to\infty}\Phi_\beta(x)=\infty$. Consequently, $\Phi_\beta$ is an example of Young function. Let $D\subseteq\R$ be a bounded (nondegenerate) interval. Recall that the \emph{(real) exponential Orlicz space} $\Osp(D;\R)=\Osp(D)=\Osp$ is defined as the linear space
\[
	\Osp(D;\R)=\left\{f\in L^0(D;\R) \left| \,\exists\lambda \in(0,\infty):\int_{D}\Phi_\beta\left(\frac{\lvert f_r\rvert}{\lambda}\right)\dd r < \infty\right.\right\}
\]
and \emph{the subspace of its finite elements} $\OspF(D;\R)=\OspF(D)=\OspF$ is defined as
\[
	\OspF(D;\R)=\left\{f\in L^0(D;\R)\left|\, \forall\lambda\in(0,\infty):\int_{D}\Phi_\beta\left(\frac{\lvert f_r\rvert}{\lambda}\right)\dd r < \infty\right.\right\}.
\]
Endowed with the \emph{Luxemburg} (or \emph{gauge}) \emph{norm}
\[
	\lVert f\rVert_{\Osp} = \inf\left\{\lambda\in(0,\infty)\,\left|\, \int_{D}\Phi_\beta\left(\frac{\lvert f_r\rvert}{\lambda}\right)\dd r \leq 1 \right.\right\}, \quad f\in \Osp,
\]
the space $\Osp$ is a Banach space with $\OspF$ being its closed subspace. We emphasize that $\OspF(D;\R) \subsetneq \Osp(D;\R)$ by, e.g., \cite[Remark 4.12.4]{PicKufJohFuc13}.

Roughly speaking, exponential Orlicz spaces measure the asymptotic growth of $L^p$-norms as $p$ increases. More precisely, there is the following equivalence (see \cite[Theorem (3.4)]{Cie93}): 
\begin{equation}
	\label{eq:OSP_equiv_norm_1}
	\norm{\,\cdot\,}{\Osp([0, 1])} \lesssim_\beta \sup_{p\in[1,\infty)}p^{-\frac{1}{\beta}} \norm{\,\cdot\,}{L^p([0, 1])} \lesssim_\beta \norm{\,\cdot\,}{\Osp([0, 1])}.
\end{equation}

It is possible to construct exponential Orlicz spaces of metric space-valued functions. Let $(E, d)$ be a~nonempty metric space. Let again $\beta\in (0,\infty)$ and let $D\subseteq\R$ be a bounded nondegenerate interval. The \emph{exponential Orlicz space of $E$-valued functions} $\Osp(D;E)$ is the linear space
\[
	\Osp(D;E)=\left\{f\in L^0(D;E)\left|\,\exists e_0\in E: d(f_\cdot, e_0)\in \Osp(D;\R)\right.\right\}
\]
and the \emph{subspace of its finite elements} $\OspF(D;E)$ is 
\[
	\OspF(D;E)=\left\{f\in L^0(D;E)\left|\,\exists e_0\in E: d(f_\cdot, e_0)\in \OspF(D;\R)\right.\right\}.
\]	
The \emph{Luxemburg norm} of $f\in\Osp(D;E)$ is defined by
\[
	\norm{f}{\Osp(D;E)}=\norm{d(f_\cdot, e_0)}{\Osp(D;\R)}
\]
for $e_0 \in E$ fixed.

There are the following observations for the space $\Osp(D;E)$ (and similarly for the space $\OspF(D;E))$:
\begin{itemize}
	\itemsep0em
	\item If $f\in\Osp(D;E)$, then, for any $e\in E$, it holds that $d(f_\cdot, e)\in\Osp(D; \R)$. In particular, as discussed in \autoref{sec:notation}, the choice of $e_0 \in E$ does not change the function space.
	\item Fix $T\in(0,\infty)$, $h\in[0, T]$, and $e\in E$. Then
	\[
		f\in \Osp([0,T];E) \quad\implies\quad d(f_{\cdot+h}, e), d(f_\cdot, f_{\cdot+h})\in \Osp([0,T-h];\R).
	\]
	\item If $E$ is complete and separable, then $(\Osp(D;E), \norm{\,\cdot\,}{\Osp(D;E)})$ is also complete.
\end{itemize}

In the rest of the article, we will only work with exponential (Besov-)Orlicz(-type) spaces and we will omit the epithet ``exponential" for simplicity. Below, we collect several properties of Orlicz spaces and Orlicz norms that will be frequently used throughout the whole text. To this end, let $T\in (0,\infty)$ be fixed. We start by a, fairly obvious, result that we often implicitly rely upon. 

\begin{lemma}
	\label{lem:luxp}
	Let $\beta\in (0,\infty)$ and $p\in (0,\infty)$. Then, for all $f\in\Osp([0,T])$, it holds
	\[
		\norm{\lvert f\rvert ^p}{\Osp([0,T])}=\normto{f}{\Ospb{p\beta}([0,T])}{p}.
	\]
\end{lemma}

The following claim is an extension of the equivalence \eqref{eq:OSP_equiv_norm_1}.

\begin{lemma}
	\label{lem:OSP_equiv_norm_2}
	Let $\beta\in (0,\infty)$. Then, for $f\in \Osp([0,T])$, it holds
	\[
		\frac{1}{(1 \vee T)}\norm{f}{\Osp([0,T])} \lesssim_{\beta} \sup_{p\in[1,\infty) }p^{-\frac{1}{\beta}}\norm{f}{L^p([0,T])} \lesssim_{\beta} (1 \vee T) \norm{f}{\Osp([0,T])}.
	\]
\end{lemma} 

\begin{proof}
	Define $f^T: [0,1]\to \R$ by $f_r^T = f_{Tr}$ for $r\in [0,1]$ if $T\geq 1$ and by $f_r^T = f_r\bm{1}_{[0,T]}(r)$ for $r\in [0,1]$ if $T<1$. Then in all the four cases ($\beta \geq 1$ and $T\geq 1$, $\beta<1$ and $T\geq 1$, $\beta\geq 1$ and $T<1$, and $\beta <1$ and $T<1$), we obtain $f^T\in\Osp([0,1])$ and, moreover, that 
	\[
		\|f^T\|_{\Osp([0,1])} \leq \|f\|_{\Osp([0,T])} \leq (1\vee T) \|f^T \|_{\Osp([0,1])}
	\]
	holds. (When $T<1$, the inequality between the norms is trivial and when $T\geq 1$, one obtains the estimate by using the fact that for every $c\geq 1$, the inequality $c(\mathrm{e}^{x}-1)\leq \mathrm{e}^{cx}-1$ holds for every $x\geq 0$.) As we also have, for every $p\in [1,\infty)$, that 
	\[ 
		\|f\|_{L^p([0,T])} = (1\vee T^\frac{1}{p})\|f^T\|_{L^p([0,1])},
	\]
	we obtain the claim by appealing to equivalence \eqref{eq:OSP_equiv_norm_1}.
\end{proof}

There is also a H\"older-type inequality for the Luxemburg norm.

\begin{lemma}
	\label{thm:orliczholder}
	Let $\beta_1,\beta_2\in(1,\infty)$ be such that $\beta_1\beta_2=\beta_1 + \beta_2$. Then the inequality
	\[
		\norm{fg}{\Ospb{1}([0,T])} \lesssim \norm{f}{\Ospb{\beta_1}([0,T])} \norm{g}{\Ospb{\beta_2}([0,T])}
	\]
	holds for all $f,g\in L^0([0,T])$. Moreover, if $\beta\in (0,\infty)$, then
	\begin{equation}
		\label{eq:orliczholder}
		\norm{fg}{\Osp([0,T])} \lesssim_\beta \normto{\lvert f\rvert^\beta}{\Ospb{\beta_1}([0,T])}{\frac{1}{\beta}} \normto{\lvert g\rvert^\beta}{\Ospb{\beta_2}([0,T])}{\frac{1}{\beta}}.
	\end{equation}
\end{lemma}

\begin{proof}
The first claim is a special case of H\"older's inequality for Orlicz spaces, see e.g.\ \cite[Theorem 7, p.\ 64]{RaoRen91}). The rest follows from \autoref{lem:luxp}.
\end{proof}

We finish this section by an embedding result from \cite[Theorem 3, p.\ 155]{RaoRen91} that can also be immediately obtained (albeit with a possibly worse constant) by repeated application of \autoref{lem:OSP_equiv_norm_2}.

\begin{proposition}
	\label{prop:osp_emb_beta}
	Let $\beta_1,\beta_2\in (0,\infty)$ be such that $\beta_1 \leq \beta_2$. If $f\in\Ospb{\beta_2}([0,T])$, then $f\in\Ospb{\beta_1}([0,T])$ and it holds
	\[
		\norm{f}{\Ospb{\beta_1}([0,T])} \lesssim_{\beta_1, \beta_2} \left(1\vee T\right) \norm{f}{\Ospb{\beta_2}([0,T])}.
	\]
\end{proposition}

\section{Exponential Besov-Orlicz-type spaces}
\label{sec:preliminaries_BO}

\subsection{Univariate Besov-Orlicz spaces}

Let us now review basic results on the Besov-Orlicz spaces of functions with values in metric spaces that will be needed for our analysis. Let us fix $T\in (0,\infty)$, $\alpha\in (0,1)$, $\beta\in (0,\infty)$,  $q\in (0,\infty]$, and a nonempty complete separable metric space $(E,d)$ for the rest of this section.

\begin{definition}
	For $f\in\Osp\onT$ and $\tau\in [0,T]$, set
	\[
		\omg{f}{\tau}=\sup_{h\in [0,\tau]} \norm{d(f_{\cdot},f_{\cdot+h})}{\Osp([0, T-h])}.
	\]
	The \emph{exponential Besov-Orlicz space} is defined as
	\[
		\BOone([0, T];E) = \left\{ f\in\Osp([0, T];E)\,\left|\,[f]_{\BOone([0,T]; E)} < \infty\right.\right\},
	\]
	where 
	\[
		[f]_{\BOone([0,T];E)} = \left\| \frac{\omg{f}{\tau}}{\tau^{\alpha}} \right\|_{\Lq}.
	\]
	We also define
	\begin{equation}
		\label{eq:1stMetric}
		d_{\BOone}(f,g) =\norm{f-g}{\Osp([0,T];E)} + [f-g]_{\BOone([0,T];E)}^{q\wedge 1}
	\end{equation}
	for $f,g\in\BOone([0,T];E)$.
\end{definition}

It can be shown that $d_{\BOone}$ is indeed a metric on $\BOone([0, T];E)$ and that $[\,\cdot\,]_{\BOone([0,T];E)}$ is a (quazi)seminorm on $\BOone([0, T];E)$ and a seminorm for $q \geq 1$. Note that since the metric space $(E,d)$ is complete and separable, the metric space $(\BOone([0,T];E), d_{\BOone([0,T];E)})$ is also complete.

Next, we summarize the embeddings between Besov-Orlicz spaces and Besov and H\"older spaces.

\begin{proposition}
	\label{thm:BOintoB}
	Let $f\in\BOone([0,T];E)$. Then for any $p\in [1,\infty)$ and $p' \in (\frac1\alpha,\infty)$, it holds
	\begin{align*}
		[f]_{\Bone([0,T];E)} &\lesssim_{\beta, p} (1 \vee T) [f]_{\BOone([0,T];E)},\\
		[f]_{C^{\alpha-1/p'}([0,T];E)} &\lesssim_{\alpha, \beta, p', q} (1 \vee T) [f]_{\BOone([0,T];E)}.		
	\end{align*}
	In particular, the inclusion $\BOone([0,T];E)\subseteq L^\infty([0, T]; E) \subseteq \OspF([0,T];E)$ holds.
\end{proposition}

\begin{proof}
	The first embedding $\BOone([0,T];E) \hook \Bone([0,T];E)$ is an immediate consequence of \autoref{lem:OSP_equiv_norm_2}. The second embedding then follows from the first one and the classical Sobolev embedding from, e.g., \cite[Theorem 3.3.1]{Tri83}. Thus, since the elements of $\BOone([0,T];E)$ are bounded, the remaining claim follows from a straightforward estimate.
\end{proof}

The Besov-Orlicz (quazi)seminorm allows a discrete characterization which will be useful below. 

\begin{proposition}
	\label{thm:equivnorm}
	For $f\in\BOone([0,T];E)$, it holds
	\begin{equation*}
		[f]_{\BOone([0,T];E)} \lesssim_{\alpha, q} \left\| \left( T2^{-n}\right)^{-\alpha} \norm{d(f_\cdot,f_{\cdot+T2^{-n}})}{\Osp([0,T-T2^{-n}])} \right\|_{\ell^q(n=1,2,\dots)}\lesssim_{\alpha, q} [f]_{\BOone([0,T];E)}.
	\end{equation*}
\end{proposition}

\begin{proof}
	The proofs of both inequalities are similar to the proofs of the corresponding inequalities in \cite[Lemma 2.2]{FriSee22}. We therefore omit the details and only discuss the nonobvious differences.
	
	The proof of the first inequality uses the continuity of the map $H: [0, T] \to [0, \infty)$ defined by $h\mapsto \norm{d(f_\cdot,f_{\cdot + h})}{\Osp([0, T-h])}$. While this continuity is almost trivial in the standard $L^p$-case, it requires a more careful approach in the case of exponential Orlicz spaces. The continuity of $H$ at $0$ follows immediately from $f\in\BOone([0,T];E)$ by contradiction. To prove its continuity in $(0, T]$, let $0\leq h_1<h_2\leq T$ be arbitrary. Then it holds that
	\begin{equation*}
		H(h_2)\leq \norm{d(f_\cdot,f_{\cdot + h_1})}{\Osp([0, T-h_2])} + \norm{d(f_{\cdot + h_1},f_{\cdot + h_2})}{\Osp([0, T-h_2])}\leq H(h_1) + H(h_2-h_1).
	\end{equation*}
	Moreover, by splitting the interval $[0,T-h_1]$ into $[0,T-h_2]$ and $[T-h_2,T-h_1]$, we deduce
	\begin{align*}
		H(h_1)&\leq \norm{d(f_\cdot,f_{\cdot + h_1})}{\Osp([0, T-h_2])} + \norm{d(f_\cdot,f_{\cdot + h_1})}{\Osp([T-h_2, T-h_1])}
		\\
		&\leq H(h_2) + \norm{d(f_{\cdot + h_1},f_{\cdot + h_2})}{\Osp([0, T-h_2])} + \norm{d(f_\cdot,f_{\cdot + h_1})}{\Osp([T-h_2, T-h_1])}
		\\
		&\leq H(h_2) + H(h_2-h_1) + \norm{d(f_\cdot,e_0)}{\Osp([T-h_2, T-h_1])} + \norm{d(f_\cdot, e_0)}{\Osp([T-(h_2-h_1), T])}.
	\end{align*}
	By combining the two estimates above, we obtain
	\begin{equation*}
		\lvert H(h_2) - H(h_1)\rvert\leq H(h_2-h_1) + \norm{d(f_\cdot,e_0)}{\Osp([T-h_2, T-h_1])} + \norm{d(f_\cdot, e_0)}{\Osp([T-(h_2-h_1), T])}.
	\end{equation*}
	The term $H(h_2-h_1)$ converges to $0$ as $(h_2-h_1)\to 0$ by the continuity of $H$ at $0$. Let us show the convergence to $0$ of the second term; the convergence of the third term follows similarly. Since $f\in\OspF([0,T];E)$ by \autoref{thm:BOintoB}, the convergence
	\begin{equation}
		\label{eq:contradiction_continuity}
		\lim_{h_2\to h_1+} \int_{T-h_2}^{T-h_1} \Phi_\beta\left(\frac{d(f_r,e_0)}{\lambda}\right)\dd{r} =0
	\end{equation}
	holds for all $\lambda\in(0,\infty)$ by the continuity of the Lebesgue integral. If $\norm{d(f_\cdot,e_0)}{\Osp([T-h_2, T-h_1])}$ does not converge to $0$, we may find $\varepsilon\in(0,\infty)$ and a sequence $\{h^n\}_{n=1}^\infty$ satisfying $h^n \searrow h_1$ such that
	\[
		\inf \left\{ \lambda \in(0,\infty) \left| \int_{T-h^n}^{T-h_1} \Phi_\beta\left(\frac{d(f_r,e_0)}{\lambda}\right)\dd{r}\right.\leq 1\right\} >\varepsilon,
	\]
	for all $n \in \N$. In particular, for all $n \in \N$, it holds
	\[
		\int_{T-h^n}^{T-h_1} \Phi_\beta\left(\frac{d(f_r,e_0)}{\frac{\varepsilon}{2}}\right)\dd{r} >1,
	\]
	which is a contradiction with \eqref{eq:contradiction_continuity}.
\end{proof}

The following proposition discusses the embeddings of Besov-Orlicz spaces.

\begin{proposition}
	\label{prop:triv_emb_BO}
	Let $\alpha_1,\alpha_2 \in (0,1)$, $\beta_1,\beta_2\in(0,\infty)$, and $q_1,q_2\in (0,\infty]$ satisfy $\alpha_1 \leq \alpha_2$, $\beta_1 \leq \beta_2$, and $q_1\geq q_2$. Then for all $f\in L^0([0,T];E)$, the following inequalities hold:
	\begin{align*}
		[f]_{\BO{\alpha_1}{\beta}{q}([0,T];E)}&\leq T^{\alpha_2-\alpha_1}[f]_{\BO{\alpha_2}{\beta}{q}([0,T];E)},
		\\
		[f]_{\BO{\alpha}{\beta_1}{q}([0,T];E)}&\lesssim_{\beta_1,\beta_2}\left(1\vee T\right)[f]_{\BO{\alpha}{\beta_2}{q}([0,T];E)},
		\\
		[f]_{\BO{\alpha}{\beta}{q_1}([0,T];E)} & \lesssim_{\alpha,q_1,q_2} [f]_{\BO{\alpha}{\beta}{q_2}([0,T];E)}.
	\end{align*}
\end{proposition}

\begin{proof}
	The first inequality follows from the definition of the Besov-Orlicz (quazi)seminorm directly. The second inequality can then be deduced from \autoref{prop:osp_emb_beta} and the third inequality can be established by appealing to \autoref{thm:equivnorm} and the monotonicity of $\ell^q$-spaces.
\end{proof}

The smoothness parameter $\alpha$ plays a more prominent role than the fine parameter $q$. In particular, a Besov-Orlicz space with smaller $\alpha$ contains Besov-Orlicz spaces of any larger $q$. More precisely, there is the following embedding result:

\begin{proposition}
	\label{cor:boeps}
	If $\varepsilon\in(0,\alpha)$, $q_1,q_2\in(0,\infty]$ satisfy $q_1\leq q_2$, and if $f\in\BO{\alpha}{\beta}{q_2}([0,T];E)$, then
	\[
		[f]_{\BO{\alpha-\varepsilon}{\beta}{q_1}([0,T];E)}\lsmT{\alpha,\varepsilon,q_1}[f]_{\BO{\alpha}{\beta}{\infty}([0,T];E)}\lesssim_{\alpha,\varepsilon,q_1,q_2}[f]_{\BO{\alpha}{\beta}{q_2}([0,T];E)}.
	\]
\end{proposition}

\begin{proof}
	We have 
	\begin{align*}
		[f]^{q_1}_{\BO{\alpha-\varepsilon}{\beta}{q_1}([0,T];E)}&\lesssim_{\alpha,\varepsilon,q_1}\sum_{n=1}^\infty \left(\left(T2^{-n}\right)^{-(\alpha-\varepsilon)} \norm{d(f_\cdot,f_{\cdot+T2^{-n}})}{\Osp([0,T-T2^{-n}])}\right)^{q_1}
		\\
		&\leq T^{\varepsilon q_1}\sup_{n\in\mathbb{N}}\left(\left(T2^{-n}\right)^{-\alpha} \norm{d(f_\cdot,f_{\cdot+T2^{-n}})}{\Osp([0,T-T2^{-n}])}\right)^{q_1}\sum_{n=1}^\infty 2^{-n\varepsilon q_1}
		\\
		&\lsmT{\alpha,\varepsilon,q_1}[f]^{q_1}_{\BO{\alpha}{\beta}{\infty}([0,T];E)}
		\\
		&\lesssim_{\alpha,q_1,q_2}[f]_{\BO{\alpha}{\beta}{q_2}([0,T];E)}^{q_1}
	\end{align*}
	where we used \autoref{thm:equivnorm} in the first and the third inequality and the last inequality follows by the third assertion of \autoref{prop:triv_emb_BO}.
\end{proof}

For $e\in E$, the subspace of $\BOone([0,T];E)$ of all functions originating from $e$ is denoted by
\begin{equation}\label{eq:affinesubspacedef}
	\BOone([0,T];E,e)=\{f\in\BOone([0,T];E)\,\left|\,f_0=e\right\}.
\end{equation}
Note that the space $\BOone([0,T];E,e)$ is well-defined since $\BOone([0,T];E) \subseteq C([0,T];E)$ by \autoref{thm:BOintoB}.
In the rest of this section, we discuss Besov-Orlicz spaces of functions with values in a separable Banach space $\V$. In particular, we establish equivalence of certain norms and recall that the seminorm $[\,\cdot\,]_{\BOone([0,T];\V)}$ is a norm on affine subspaces of $\BOone([0,T];\V)$ of the form \eqref{eq:affinesubspacedef} for $q \geq 1$. This will be important for controlled rough paths in subsequent sections.

\begin{proposition}
	\label{lem:affinesub}
	Let $(\V,\lvert\,\cdot\,\rvert_\V)$ be a separable Banach space.
	\begin{enumerate}[label=\roman*)]
		\item The mappings
		\begin{align*}
			N_{\Phi_{\beta}}: \BOone([0,T];\V)\to\R: f&\mapsto(\norm{f}{\Osp([0,T];\V)}+[f]_{\BOone([0,T];\V)}),
			\\
			N_0: \BOone([0,T];\V)\to\R: f&\mapsto(\lvert f_0\rvert_\V + [f]_{\BOone([0,T];\V)}),
			\\
			N_\infty: \BOone([0,T];\V)\to\R: f&\mapsto(\norm{f}{L^\infty([0,T];\V)} + [f]_{\BOone([0,T];\V)})
		\end{align*}
		are equivalent quazinorms on $\BOone([0,T];\V)$, resp.\ equivalent norms if $q \geq 1$, and, for $f \in \BOone([0,T];\V)$, it holds
		\begin{align*}
			TN_0(f)\lsmT{\alpha, \beta, q}N_{\Phi_{\beta}}(f)&\lsmT{\alpha, \beta, q}N_0(f),\\
			TN_\infty(f)\lsmT{\alpha, \beta, q}N_{\Phi_{\beta}}(f)&\lsmT{\alpha, \beta, q}N_\infty(f).
		\end{align*}
		Moreover, for all $p\in(\frac1\alpha,\infty)$, it holds
		\[
			\norm{f}{L^\infty([0,T];\V)} \lesssim_{\alpha,\beta,q,p} \lvert f_0\rvert_{\V} + T^{\alpha-\frac1p}(1\vee T)[f]_{\BOone([0,T];\V)}.
		\]
		\item Let $v\in\V$ 
		and set
		\[
			\rho(f,g)=[f-g]^{q\wedge 1}_{\BOone([0,T];\V)}
		\]
		for $f,g\in \BOone([0,T];\V,v)$. Then $\rho$ is a metric on $\BOone([0,T];\V,v)$ and, $[\,\cdot\,]_{\BOone([0,T];\V)}$ is a quazinorm, resp.\ norm if $q\geq 1$, on $\BOone([0,T];\V,v)$. Moreover, $\BOone([0,T];\V,v)$ is a closed (affine) subspace of $\BOone([0,T];\V)$.
	\end{enumerate}
\end{proposition}

It is important to note that the constants in \autoref{lem:affinesub} depend on $T$ in a nondecreasing manner.

\begin{proof}
	The proof of the equivalence of the norms mirrors the standard proofs of the equivalence of the norms on spaces of H\"older-continuous functions and it is therefore omitted. The remaining estimate in i) can be established by
	\[
		\|f\|_{L^\infty([0,T];\V)} \leq |f_0|_\V+ \sup_{r\in [0,T]} |f_r-f_0|_\V \leq |f_0|_\V+ T^{\alpha-\frac{1}{p}} [f]_{C^{\alpha-1/p}([0,T];\V)}
	\]
	and the embedding into H\"older-continuous functions from \autoref{thm:BOintoB}. The properties of $\rho$ and $f\mapsto[f]_{\BOone([0,T];\V)}$ follow from the definitions and the restriction $f_0 = v$.
\end{proof}

\subsection{Multivariate Besov-Orlicz-type spaces}

In this section, Besov-Orlicz-type spaces of multivariate maps suitable for rough path analysis are defined and some of their properties are collected. In the whole section, let us fix a nonempty separable Banach space $(\V,|\cdot|_\V)$, $d\in \{2,3\}$, $T\in (0,\infty)$, $\alpha\in (0,\infty)$, $\beta\in (0,\infty)$, and $q\in (0,\infty]$.

\begin{definition}
	For $\varXi\in L^{0;d}([0,T];\V)$ and $\tau\in[0, T]$, set
	\begin{equation}
		\omgd{\varXi}{\tau}{d} = \begin{cases}
			\displaystyle \sup_{h\in[0, \tau]}\norm{\absh{\varXi_{r, r+h}}{\V}}{\Osp([0, T-h])},& \quad d=2,\\
			\displaystyle \sup_{\theta\in[0, 1]}\sup_{h\in[0, \tau]}\norm{\absh{\varXi_{r, r+\theta h, r+h}}{\V}}{\Osp([0, T-h])},& \quad d=3.
		\end{cases}
	\end{equation}
	We define the \emph{(multivariate) exponential Besov-Orlicz-type space} by
	\[
		\BOd{d}([0, T];\V) = \left\{\varXi\in L^{0;d}([0,T];\V)\left|\,\norm{\varXi}{\BOd{d}([0,T];\V)} < \infty\right.\right\},
	\]
	where
	\[
		\norm{\varXi}{\BOd{d}([0,T];\V)} = \left\| \frac{\omgd{\varXi}{\tau}{d}}{\tau^{\alpha}} \right\|_{\Lq}.
	\]
	For $\varXi, \tilde{\varXi} \in \BOd{d}([0,T];\V)$, we also define
	\begin{equation}\label{eq:2ndMetric}
		d_{\BOd{d}([0,T];\V)}(\varXi, \tilde{\varXi}) =\normto{\varXi-\tilde{\varXi}}{\BOd{d}([0,T];\V)}{q\wedge 1}.
	\end{equation}
\end{definition}

One can show that $\norm{\,\cdot\,}{\BOd{d}([0,T];\V)}$ is a quazinorm, resp.\ norm if $q \geq 1$, on $\BOd{d}([0,T];\V)$ and that $d_{\BOd{d}([0,T];\V)}$ is a metric on $\BOd{d}([0, T];\V)$.  Note that if $\V$ is a separable Banach space, then the space $\BOd{d}([0,T];\V)$ is also Banach for both $d \in \{ 2, 3 \}$.

First, we summarize the relation of the Besov-Orlicz-type spaces to the Besov-type and H\"older-type spaces defined in \autoref{sec:preliminaries}.

\begin{proposition}
	\label{thm:emb_BO_B_two_param}
	Let $\varXi \in\BOd{2}([0,T];\V)$, then for all $p \in [1, \infty)$ it holds
	\[
		\norm{\varXi}{\Besov{\alpha;2}{p}{q}([0,T];\V)}\lesssim_{\beta,p}(1\vee T)\norm{\varXi}{\BOd{2}([0,T];\V)}.
	\]
	If $\V = \R^m$ and if, additionally, $\varXi$ satisfies
	\begin{equation}
		\label{eq:deltaMzeta}
		\lvert\delta \varXi_{s,u,t}\rvert_{\R^m}\leq M\left(\left(\left(u-s\right)\wedge\left(t-u\right)\right)^\theta\left(\left(u-s\right)\vee\left(t-u\right)\right)^{1-\theta}\right)^{\alpha-\frac1p}
	\end{equation}
	for some $M\in(0,\infty)$, $\theta \in (0, \frac12]$, $p \in (\frac1\alpha, \infty)$, and all $(s,u,t)\in\triangle^3[0, T]$, then $\varXi$ has a continuous version and it holds
	\[	
		\norm{\varXi}{C^{\alpha-1/p;2}([0,T];\R^m)} \lesssim_{\theta,\alpha,\beta,p} (1\vee T) \norm{\varXi}{\BOd{2}([0,T];\R^m)}+M.
	\]
\end{proposition}

\begin{proof}
	The first claim follows from the definitions of the respective norms and \autoref{lem:OSP_equiv_norm_2} immediately. The second claim can be established by \cite[Proposition 2.7]{FriSee22} and the first claim.
\end{proof}

Next, we give an interpolation result whose proof is similar to that of \cite[Lemma 2.7]{FriSee22} and is therefore omitted.

\begin{lemma}\label{lem:interpolation}
	Let $\gamma\in (\alpha,\infty)$, $\beta^\prime\in(\beta,\infty)$, and $\delta\in (0,\infty)$ satisfy $\delta(1-\frac{\beta}{\beta'}) + \frac{\gamma\beta}{\beta'}-\alpha>0$.
	Then
	\[
		\norm{\varXi}{\BO{\alpha;2}{\beta^\prime}{q}([0,T];\V)} \lesssim_{\alpha, \gamma, \beta, \beta', q, \delta} T^{\delta(1-\frac{\beta}{\beta'})+\frac{\gamma\beta}{\beta'}-\alpha} \normto{\varXi}{C^{\delta;2}([0,T];\V)}{1-\frac{\beta}{\beta'}}\normto{\varXi}{\BO{\gamma;2}{\beta}{q}([0,T];\V)}{\frac{\beta}{\beta'}}
	\]
	holds for any $f\in L^{0;d}([0,T];\V)$.
\end{lemma}

In the rest of the section, we discuss embeddings of the Besov-Orlicz-type spaces defined above. The monotonicity of Besov-Orlicz-type spaces with respect to the smoothness parameter $\alpha$ and the integrability parameter $\beta$ can be established analogously as in \autoref{prop:triv_emb_BO} and therefore, we skip the proof.

\begin{proposition}
	\label{prop:trivemb}
	If $\alpha_1,\alpha_2,\beta_1,\beta_2\in (0,\infty)$ satisfy $\alpha_1 \leq \alpha_2$ and $\beta_1 \leq \beta_2$, then 
	\begin{align*}
		\norm{\varXi}{\BO{\alpha_1;d}{\beta}{q}([0,T];\V)}&\leq T^{\alpha_2-\alpha_1}\norm{\varXi}{\BO{\alpha_2;d}{\beta}{q}([0,T];\V)},
		\\
		\norm{\varXi}{\BO{\alpha;d}{\beta_1}{q}([0,T];\V)}&\lesssim_{\beta_1,\beta_2} (1\vee T)\norm{\varXi}{\BO{\alpha;d}{\beta_2}{q}([0,T];\V)}
	\end{align*}
	hold for $\varXi \in L^{0;d}([0,T];\V)$. 
\end{proposition}

Compared to the standard Besov-Orlicz, resp.\ Besov, spaces the definition of the Besov-Orlicz-type, resp.\ Besov-type, space does not enforce the boundedness of $\omgd{\varXi}{\cdot}{d}$ on $[0, T]$. Thus, the monotonicity of the Besov-Orlicz-type spaces with respect to the fine parameter $q$ does not hold without additional assumptions. We demonstrate this in the following example. 

\begin{example}
	Let $\theta\in(1,\infty)$ and consider $\varXi^\theta:\triangle^2[0,T]\to\R$ such that
	\[
		\varXi^\theta_{r,r+h} = \frac{(T-h)^{-\theta} \bm{1}_{[\frac{T}{2},T]}(h)}{\norm{1}{\Osp([0, T-h])}}, \quad h\in [0,T),\ r\in [0,T-h].
	\]
	Then $\varXi^\theta \in \BOd{2}([0,T];\R)$ if and only if $q<\frac{1}{\theta-1}$. Hence, choosing $0<q_2<\frac{1}{\theta-1}<q_1$, we observe $\|\varXi^\theta\|_{\BO{\alpha;2}{\beta}{q_1}([0,T];\R)}=\infty$ while $\|\varXi^\theta\|_{\BO{\alpha;2}{\beta}{q_2}([0,T];\R)}<\infty$.
\end{example}

\begin{proposition}
	\label{prop:emb_BOq1q2}  
	Let $q_1,q_2\in (0,\infty]$ be such that $q_2 \leq q_1$ and let $\varXi \in \BO{\alpha;2}{\beta}{q_2}([0,T];\R^m)$ be such that $\omgd{\varXi}{\cdot}{d}$ is bounded on $[0, T]$. Then $\varXi \in \BO{\alpha;2}{\beta}{q_1}([0,T];\R^m)$.

	If, in addition, \eqref{eq:deltaMzeta} holds for some $M\in(0,\infty)$, $\theta \in (0, \frac12]$, $p \in (\frac1\alpha, \infty)$, and all $(s,u,t)\in\triangle^3[0, T]$, then there is the following estimate:
	\[
		\norm{\varXi}{\BO{\alpha;2}{\beta}{q_1}([0,T];\R^m)}\lesssim_{\alpha,\beta,p,q_1,q_2}^T\norm{\varXi}{\BO{\alpha;2}{\beta}{q_2}([0,T];\R^m)}+MT^{-\frac1p}.
	\]
\end{proposition}

With a slight abuse of terminology, we will often refer to the result of \autoref{prop:emb_BOq1q2} as ``the inclusion $\BO{\alpha;2}{\beta}{q_2} \subseteq \BO{\alpha;2}{\beta}{q_1}$" and we will not mention the additional assumptions explicitly.

\begin{proof}
	Let $q_2<\infty$. Since $\varXi \in \BO{\alpha;2}{\beta}{q_2}([0,T];\R^m)$, the function $\tau \mapsto \omgd{\varXi}{\tau}{d}\tau^{-\alpha}$ is bounded near $0$ and thus, by the assumption on boundedness of $\omgd{\varXi}{\cdot}{d}$, we obtain $\varXi \in \BO{\alpha;2}{\beta}{\infty}([0,T];\R^m)$. For $q_1 < \infty$, the inclusion $\BO{\alpha;2}{\beta}{q_2}([0,T];\R^m) \subseteq \BO{\alpha;2}{\beta}{q_1}([0,T];\R^m)$ then follows by the straightforward estimate
	\[
		\norm{\varXi}{\BO{\alpha;2}{\beta}{q_1}([0,T];\R^m)} \leq \norm{\varXi}{\BO{\alpha;2}{\beta}{\infty}([0,T];\R^m)}^{1-q_2/q_1} \norm{\varXi}{\BO{\alpha;2}{\beta}{q_2}([0,T];\R^m)}^{q_2/q_1} < \infty.	
	\]
	
	Let us show the explicit bound for $q_1 < \infty$; the case $q_1 = \infty$ can be established similarly. First, from \eqref{eq:deltaMzeta}, we deduce that
	\[
		\|\varXi_{\cdot,\cdot +h}\|_{\Osp([0,T-h])} \leq \| \varXi_{\cdot,\cdot + \frac{h}{2}}\|_{\Osp([0,T-h])} + \|\varXi_{\cdot + \frac{h}{2},\cdot + h}\|_{\Osp([0,T-h])} + MT^{\alpha-\frac{1}{p}} \|1\|_{\Osp([0,T-h])}
	\]
	holds for every $h\in [0,T]$. Since $\|1\|_{\Osp([0,T-h])} \leq C_\beta (1\vee T)$ for some $C_\beta\in(0,\infty)$ by a direct calculation, the above implies
	\begin{equation}
		\label{eq:omega_at_T_bound}
		\omega_{\Phi_\beta}^d(\varXi,T) \leq 2\omega_{\Phi_\beta}^d\left(\varXi,\frac{T}{2}\right) + M_{\beta, T} T^{\alpha-\frac{1}{p}},
	\end{equation}
	where we write $M_{\beta, T} = M C_\beta (1\vee T)$ for brevity.
	
	Next, we estimate $\|\varXi\|_{\BO{\alpha;d}{\beta}{q_1}([0,T];\R^m)}$ by a discrete sum similarly as in \autoref{thm:equivnorm}. In particular, we split the interval $[0,T]$ in the outer integral in the definition of $\|\varXi\|_{\BO{\alpha;d}{\beta}{q_1}([0,T];\R^m)}$ into intervals $[T2^{-n-1},T2^{-n}]$, $n\in \mathbb{N}_0$, and use the estimates 
	\[
		\tau^{-\alpha q_1} \leq (T2^{-n-1})^{-\alpha q_1}\quad \mbox{and} \quad \omega_{\Phi_\beta}^d(\varXi,\tau) \leq \omega_{\Phi_\beta}^d(\varXi,T2^{-n})
	\]
	that holds for $\tau\in [T2^{-n-1},T2^{-n}]$ to obtain
	\begin{equation}
		\label{eq:upper_BOd}
		\|\varXi\|_{\BO{\alpha;d}{\beta}{q_1}([0,T];\R^m)}	\leq (\log 2)^\frac{1}{q_1}\left(\frac{2}{T}\right)^\alpha \left( \sum_{n=0}^\infty \left( 2^{\alpha n}\omega_{\Phi_\beta}^d(\varXi,T2^{-n})\right)^{q_1}\right)^\frac{1}{q_1}.
	\end{equation}
	Applying \eqref{eq:omega_at_T_bound} to the zeroth term of the sum in \eqref{eq:upper_BOd} and recalling the monotonicity of $\ell^q$-spaces and the inequality $(a+b)^q \leq 2^q(a^q + b^q)$ that holds for $a,b\in[0,\infty)$ and $q\in[0,\infty)$, one easily derives
	\[
		\|\varXi\|_{\BO{\alpha;d}{\beta}{q_1}([0,T];\R^m)} 
		\lesssim_{\alpha,q_1,q_2} T^{-\alpha} \left( \sum_{n=1}^\infty \left(2^{\alpha n}\omega_{\Phi_\beta}^d(\varXi,T2^{-n})\right)^{q_2} \right)^\frac{1}{q_2}+ M_{\beta,T} T^{-\frac{1}{p}}.
	\]
	The desired claim now follows from the inequality 
	\[
		\sum_{n=1}^\infty \left(2^{\alpha n}\omega_{\Phi_\beta}^d(\varXi,T2^{-n})\right)^{q_2} \leq \frac{(2T)^{\alpha q_2}}{\log 2} \|\varXi\|_{\BO{\alpha;d}{\beta}{q_2}([0,T];\R^m)}^{q_2},
	\]
	which can be proved similarly as \eqref{eq:upper_BOd}.
\end{proof}

\section{Exponential Besov-Orlicz rough paths}
\label{sec:rps}

In this section, we study the exponential Besov-Orlicz rough paths and the related results such as a sewing lemma, rough integral and its stability. Our treatment is not essentially different from the ones in, e.g., Friz and Hairer \cite{FriHai20} and Friz and Seeger \cite{FriSee22}. However, in contrast to the latter reference, we consider only compositions of controlled rough paths with $C^3_b$ functions, see the discussion above \autoref{thm:c2comp} below.

Let $\alpha\in(\frac{1}{3},\frac{1}{2}]$, $\beta\in(0,\infty)$, and $q\in(0,\infty]$ be fixed in the whole section.

\subsection{The space of exponential Besov-Orlicz rough paths}

We begin by the definition of the Besov-Orlicz rough paths and provide some basic estimates.

\begin{definition}
	Let $X\in\BOone([0,T];\Rborpcodo)$ and $\Xbb\in\BOtwo([0,T];\R^{n\times n})$ satisfy the so-called \emph{Chen's relations}, i.e.\
	\[
		\mathbb{X}_{s, t} - \mathbb{X}_{s, u} - \mathbb{X}_{u, t} = (\delta X_{s,u})(\delta X_{u,t})^\top,\quad (s,u,t) \in \triangle^3[0, T].
	\] Then the pair $\X=(X,\Xbb)$ is called an \emph{exponential Besov-Orlicz rough path (over $\Rborpcodo$)}, or an \emph{exponential Besov-Orlicz rough path lift of $X$}. We denote the space of all exponential Besov-Orlicz rough paths by $\BOrp=\BOrp([0,T];\Rborpcodo)$.
\end{definition}

In what follows, we show that a Besov-Orlicz rough path can be identified with a function that takes values in a truncated tensor algebra of level 2 of Besov-Orlicz regularity. Recall that the truncated tensor algebras of level 1 and of level 2 are defined by 
\begin{equation*}
	T_1^{(1)}(\Rborpcodo) =  \{(1,b)\,|\, b\in\Rborpcodo\}\quad\text{and}\quad T_1^{(2)}(\Rborpcodo) =  \{(1,b,c)\,|\, b\in\Rborpcodo, c\in\R^{n\times n}\},
\end{equation*}
respectively. We also define
\begin{align*}
	\delta_\lambda&:T_1^{(2)}(\Rborpcodo)\to T_1^{(2)}(\Rborpcodo):(1,b,c)\mapsto(1,\lambda b,\lambda^2 c),
	\\
	\tau_1&:T_1^{(2)}(\Rborpcodo)\to T_1^{(1)}(\Rborpcodo):(1,b,c)\mapsto(1,b),
	\\
	\tau_2&:T_1^{(2)}(\Rborpcodo)\to T_1^{(2)}(\Rborpcodo):(1,b,c)\mapsto(1,b,c).
\end{align*}
and recall that the tensor product $\otimes$ on $T_1^{(2)}(\Rborpcodo)$ is defined by 
\[
    (1,b,c)\otimes(1,\tilde{b},\tilde{c})=(1,b+\tilde{b},c+\tilde{c}+b\tilde{b}^\top)
\]
for $(1,b,c),(1,\tilde{b}, \tilde{c})\in T_1^{(2)}(\Rborpcodo)$. Clearly, $(1,0,0)$ is the unit element with respect to $\otimes$ and $(1,b,c)^{-1}=(1,-b,-c+bb^\top)$. Moreover, we define
\begin{equation*}
	N:T_1^{(2)}(\Rborpcodo)\to\R:(1,b,c)\mapsto\max\left\{\absh{b}{\Rborpcodo},\sqrt{2\absh{c}{\R^{n\times n}}}\right\}
\end{equation*}
and
\begin{equation*}
	\vertiii{\cdot}{}{}:T_1^{(2)}(\Rborpcodo)\to\R:\mathbf{x}\mapsto\frac{1}{2}(N(\mathbf{x})+N(\mathbf{x}^{-1})).
\end{equation*}
It is possible to show that the map
\[
	d_{T_1^{(2)}(\Rborpcodo)}:T_1^{(2)}(\Rborpcodo)\times T_1^{(2)}(\Rborpcodo)\to\R:(\mathbf{x,\mathbf{y}})\mapsto\vertiii{\mathbf{x}^{-1}\otimes\mathbf{y}}{}{}
\]
is a metric on the space $T_1^{(2)}(\Rborpcodo)$ and that this space when endowed with metric $d_{T_1^{(2)}(\Rborpcodo)}$ is  separable.

On the other hand, on the space $\BOrp$, we define the map
\[
	\vertiii{\X}{\BOrp}{} = [X]_{\BOone} + \lVert \mathbb{X}\rVert_{\BOtwo}^{\frac{1}{2}}
\]
for $\X=(X,\Xbb)\in\BOrp$. With a slight abuse of terminology, we call $\vertiii{\cdot}{\BOrp}{}$ \emph{the exponential Besov-Orlicz rough path norm} even though it is homogeneous only in the sense that
\[
	\vertiii{\delta_\lambda\X}{\BOrp}{}=\lvert\lambda\rvert\vertiii{\X}{\BOrp}{}
\]
holds for any $\lambda\in\R$ and $\X\in\BOrp$. For $\X = (X,\Xbb), \mathbf{\tilde{X}}=(\tilde{X},\tilde{\Xbb})\in\BOrp$, we set
\[
	\rho_{\BOrp}(\X, \mathbf{\tilde{X}}) = d_{\BOone}(X, \tilde{X}) + d_{\BOtwo}(\mathbb{X}, \tilde{\mathbb{X}}),
\]
where the metrics $d_{\BOone}$ and $d_{\BOtwo}$ are defined by \eqref{eq:1stMetric} and \eqref{eq:2ndMetric}, respectively. We call $\rho_{\BOrp}$ \emph{the exponential Besov-Orlicz rough path metric}.

Finally, let us note that, for $\X = (X,\Xbb) \in \BOrp$, the notation $\X^{(1)} := \delta X$ and $\X^{(2)} := \Xbb$ will sometimes be used.

\begin{lemma}\label{lem:deltaXdeltaXdiffleqM}
	Let $\X = (X,\Xbb)\in\BOrp$ and $p \in (\frac1\alpha, \infty)$. Then, for all $(s,u,t)\in\triangle^3[0,T]$, it holds
	\[
		\lvert\delta\Xbb_{s,u,t}\rvert_{\R^{n\times n}}\lsmT{\alpha,\beta,p,q}[X]_{\BOone}^2\left(((u-s)\wedge (t-u))^\frac{1}{2}((u-s)\vee(t-u))^\frac{1}{2}\right)^{2(\alpha-\frac{1}{p})}
	\]
	and
	\begin{multline*}
		\lvert\delta(\Xbb-\tilde{\Xbb})_{s,u,t}\rvert_{\R^{n\times n}} \lsmT{\alpha,\beta,p,q}\left([X]_{\BOone}[X-\tilde{X}]_{\BOone} + [X-\tilde{X}]_{\BOone}[\tilde{X}]_{\BOone}\right)
		\\
		\cdot\left(((u-s)\wedge (t-u))^\frac{1}{2}((u-s)\vee(t-u))^\frac{1}{2}\right)^{2(\alpha-\frac{1}{p})}.
	\end{multline*}
\end{lemma}

\begin{proof}
	The first claim follows from Chen's relations, the compatibility of the $\R^{n\times n}$-norm and the embedding $\BOone \hook C^{\alpha-1/p}$ from \autoref{thm:BOintoB} by
	\begin{align*}
		\lvert\delta\Xbb_{s,u,t}\rvert_{\R^{n\times n}}&\leq[X]_{C^{\alpha-1/p}}^2\lvert u-s\rvert^{\alpha-\frac{1}{p}}\lvert t-u\rvert^{\alpha-\frac{1}{p}}
		\\
		&\lsmT{\alpha,\beta,p,q}[X]_{\BOone}^2\left(((u-s)\wedge (t-u)^\frac{1}{2}((u-s)\vee(t-u))^\frac{1}{2}\right)^{2(\alpha-\frac{1}{p})}.
	\end{align*}
	Similarly, we deduce
	\[
		\delta(\Xbb-\tilde{\Xbb})_{s,u,t} = (\delta X_{s,u})(\delta(X-\tilde{X})_{u,t})^\top + (\delta(X-\tilde{X})_{s,u})(\delta\tilde{X}_{u,t})^\top
	\]
	and establish the second claim analogously.
\end{proof}

\begin{lemma}\label{lem:preHolderRPEmbs}
	Let $\X,\mathbf{\tilde{X}}\in\BOrp$ and $p \in (\frac1\alpha, \infty)$. Then, for $k=1,2$, it holds
	\begin{align*}
		\normto{\X^{(k)}}{C^{k(\alpha-1/p);2}}{\frac{1}{k}} &\lsmT{\alpha,\beta,p,q}\vertiii{\tau_k\X}{\BOrp}{},
		\\
		\norm{\X^{(k)}-\mathbf{\tilde{X}}^{(k)}}{C^{k(\alpha-1/p);2}} &\lsmT{\alpha,\beta,p,q}\sum_{j=1}^k\left(\vertiii{\tau_k\X}{\BOrp}{}\vee\vertiii{\tau_k\mathbf{\tilde{X}}}{\BOrp}{}\right)^{k-j}\norm{\X^{(j)} - \mathbf{\tilde{X}}^{(j)}}{\BO{j\alpha;2}{\beta/j}{q/j}}.
	\end{align*}
\end{lemma}

\begin{proof}
	The result is an immediate consequence of similar bounds for Besov rough paths from \cite[Proposition 5.1]{FriSee22} thanks to the embeddings $\BOone \hook B^\alpha_{p,q}$ from \autoref{thm:BOintoB} and $\BO{j\alpha;2}{\beta/j}{q/j} \hook \Besov{j\alpha;2}{p/j}{q/j}$ from \autoref{thm:emb_BO_B_two_param}.
\end{proof}

\begin{lemma}\label{lem:HolderRPEmbs}
	Let $\X,\mathbf{\tilde{X}}\in\BOrp$ and $p \in (\frac1\alpha, \infty)$. Then it holds
	\begin{align*}
		\norm{\Xbb}{\BOd{2}} &\lsmT{\alpha,\beta,p,q}(T^{\alpha-\frac{1}{p}}\vee T^{\alpha-\frac{2}{p}})\vertiii{\X}{\BOrp}{2},
		\\
		\norm{\Xbb-\tilde{\Xbb}}{\BOd{2}} &\lsmT{\alpha,\beta,p,q}(T^{\alpha-\frac{1}{p}}\vee T^{\alpha-\frac{3}{p}})\sum_{j=1}^2\left(\vertiii{\X}{\BOrp}{}\vee\vertiii{\mathbf{\tilde{X}}}{\BOrp}{}\right)^{2-j}\norm{\X^{(j)} - \mathbf{\tilde{X}}^{(j)}}{\BO{j\alpha;2}{\beta/j}{q/j}}.
	\end{align*}
\end{lemma}

\begin{proof}
	The proof is analogous to \cite[Lemma 5.1]{FriSee22} with the appropriate Besov-Orlicz variants of the interpolation result from \autoref{lem:interpolation} and the H\"older-type bound \autoref{lem:preHolderRPEmbs} instead of \cite[Lemma 2.7]{FriSee22} and \cite[Proposition 5.1]{FriSee22}, respectively. We establish the first claim in detail to illustrate the origin of the negative powers of $T$.
	
	From \autoref{lem:interpolation} (with $\delta=2\alpha-\frac{2}{p}$, $\gamma=2\alpha$, $\frac{\beta}{2}$ in place of $\beta$, and $\beta'=\beta$), the inclusion $\BO{\alpha;2}{\beta}{q_2} \subseteq \BO{\alpha;2}{\beta}{q_1}$ for $q_2 \leq q_1$ from \autoref{prop:emb_BOq1q2} (with $M = [X]^2_{\BOone}$ from \autoref{lem:deltaXdeltaXdiffleqM}), and from the H\"older bound in \autoref{lem:preHolderRPEmbs}, we obtain
	\begin{align*}
			\norm{\Xbb}{\BOd{2}} &\lesssim_{\alpha,\beta,p,q} T^{\alpha-\frac{1}{p}}\normto{\Xbb}{C^{2\alpha-2/p;2}}{\frac{1}{2}}\normto{\Xbb}{\BO{2\alpha;2}{\beta/2}{q}}{\frac{1}{2}}\\
		&\lsmT{\alpha,\beta,p,q} T^{\alpha-\frac{1}{p}} \vertiii{\X}{\BOrp}{} \left(\norm{\Xbb}{\BOtwo}+[X]^2_{\BOone}T^{-\frac{2}{p}}\right)^{\frac{1}{2}}.
	\end{align*}
	The claim then follows from the definition of the norm on $\BOrp$.
\end{proof}

The result that connects the space of Besov-Orlicz rough paths to the space of $T_1^{(2)}(\Rborpcodo)$-valued functions of exponential Besov-Orlicz regularity is given now.

\begin{proposition}
	The following claims hold:
	\begin{enumerate}[label=\roman*)]
		\item If $(X,\Xbb)\in\BOrp([0,T];\R^n)$, then the $T_1^{(2)}(\Rborpcodo)$-valued path $\mathbf{x}_\cdot = (1, \delta X_{0,\cdot}, \Xbb_{0,\cdot})$ satisfies $\mathbf{x}\in\BOone([0, T]; T_1^{(2)}(\Rborpcodo), (1,0,0))$ and
		\begin{equation}\label{eq:equivalenceRPMP}
		\vertiii{\X}{\BOrp([0,T];\R^n)}{}\asymp_{q} [\mathbf{x}]_{\BOone([0, T]; T_1^{(2)}(\Rborpcodo))}.
		\end{equation}
		\item Conversely, if $\mathbf{x}\in\BOone([0, T]; T_1^{(2)}(\Rborpcodo), (1,0,0))$, define $\X = (X,\Xbb)$ by $(1, \delta X_{s,t}, \Xbb_{s,t})=\mathbf{x}_s^{-1}\otimes\mathbf{x}_t$ for all $(s,t)\in\triangle^2[0,T]$. Then $\X\in\BOrp([0,T];\R^n)$ and \eqref{eq:equivalenceRPMP} holds.
	\end{enumerate}
\end{proposition}

\begin{proof}
In the first claim, equivalence \eqref{eq:equivalenceRPMP} can be established from the definitions above by elementary estimates and \autoref{lem:luxp}. From \autoref{lem:preHolderRPEmbs}, we deduce that $d_{T_1^{(2)}(\Rborpcodo)}((1,\delta X_{0,\cdot},\Xbb_{0,\cdot}),(1,0,0))$ is bounded on $[0, T]$ and thus $\mathbf{x}\in\Osp([0,T];T_1^{(2)}(\Rborpcodo))$ by \autoref{lem:affinesub}.

Regarding the second claim, note that $\X$ (and in particular $X$) is indeed well-defined since, if $\mathbf{x}_\cdot = (1, b_\cdot, c_\cdot)$ for some $b: [0,T] \to \Rborpcodo$ and $c: [0,T] \to \R^{n\times n}$, then
\[
	\mathbf{x}_s^{-1}\otimes\mathbf{x}_t = (1, \delta b_{s,t}, \delta c_{s,t}-b_s(\delta b_{s,t})^\top), \quad (s,t) \in \triangle^2[0,T],
\]
hence $X=b$. The proof of the second claim follows similarly as the proof of the first one.
\end{proof}

\subsection{Brownian motion as an exponential Besov-Orlicz rough path and other examples}

Below, we show that any path $X\in\BOone([0,T];\R)$ can be lifted to a Besov-Orlicz rough path 
	\[
		\X=(X,\Xbb)\in\BOrp([0,T];\R).
	\]	

\begin{theorem}\label{thm:lift}
	Let $X\in\BOone([0,T];\R)$ and let $F\in\BOtwo([0,T];\R)$ be an additive function, i.e.\ let $\delta F_{s,u,t}=0$ hold for all $(s,u,t) \in \triangle^3[0,T]$. Define $\Xbb: \triangle^2[0,T] \to \R$ by
	\[
		\Xbb_{s,t}=\frac{(\delta X_{s,t})^2}{2} + F_{s,t}, \quad (s,t)\in\triangle^2[0,T].
	\]
	Then $\Xbb$ belongs to the space $\BOtwo([0,T];\R)$ and satisfies Chen's relations.
\end{theorem}

\begin{proof}
	Since $\norm{(\delta X)^2}{\BOtwo} = [X]^2_{\BOone}$ holds by the definitions of the respective quazi(semi)norms and \autoref{lem:luxp}, the desired bound can be established by the (quazi-)triangle inequality as
	\[
		\norm{\Xbb}{\BOtwo} \lesssim_q \norm{(\delta X)^2}{\BOtwo} + \norm{F}{\BOtwo} = [X]^2_{\BOone} + \norm{F}{\BOtwo} < \infty.
	\]
	Chen's relations can be verified directly.
\end{proof}

Let us fix a probability basis $(\Omega, \mathscr{F}, (\mathscr{F}_t)_{t\in[0,1]}, \mathsf{P})$ rich enough to carry any of the processes below.
\autoref{thm:lift} enables us to construct a rough path lift for any ($\R$-valued) continuous local martingale $X$ with Lipschitz continuous quadratic variation. It is well-known that for such a process, it holds that $X(\omega)\in\BO{1/2}{2}{\infty}([0,1];\R)$ for a.a.\ $\omega\in\Omega$; see \cite[Theorem 4.1]{OndSimKup18}. Set
\begin{equation*}
	\Xbbito_{s,t} = \frac{(\delta X_{s,t})^2}{2} - \frac{\delta \langle X\rangle_{s,t}}{2}\quad\mathrm{and}\quad\Xbbstrat_{s,t} = \frac{(\delta X_{s,t})^2}{2}
\end{equation*}
for $(s,t)\in\triangle^2[0,1]$ where $\langle X\rangle$ stands for the (probabilistic) quadratic variation of $X$. The mapping~$\Xbbito_{s,t}$ coincides with the Itô integral $\int_{s}^{t} (\delta X_{s,r}) \dd X_r$ while the mapping $\Xbbstrat_{s,t}$ coincides with the Stratonovich integral $\int_{s}^{t} (\delta X_{s,r})\circ\dd X_r$. By appealing to \autoref{thm:lift} with $F^{\mathrm{It\hat{o}}}_{s,t}=-\tfrac12 \delta \langle X\rangle_{s,t}$ and $F^{\mathrm{Strat}}_{s,t}=0$, we obtain
\begin{align*}
	\Xito(\omega)&=(X(\omega),\Xbbito(\omega))\in\mathscr{B}^{1/2}_{\Phi_2, \infty}\quad\text{for a.a.\ } \omega\in\Omega,
	\\
	\Xstrat(\omega)&=(X(\omega),\Xbbstrat(\omega))\in\mathscr{B}^{1/2}_{\Phi_2, \infty}\quad\text{for a.a.\ } \omega\in\Omega.
\end{align*}
A canonical example of these constructions is the standard scalar Brownian motion $W$ for which we have that $W(\omega)\in B^{1/2}_{\Phi_2, \infty}([0,1];\R)$ holds for a.a. $\omega\in\Omega$ (see \cite[Theorem 5.8]{Cie93} for a stand-alone proof).
As far as Gaussian processes with lower regularity are concerned, the scalar fractional Brownian motion $W^H$ with Hurst parameter $H\in(\frac{1}{3},\frac{1}{2}]$ can be mentioned. By \cite[Corollary 5.3]{Ver09} (see also \cite{CieKerRoy93}), we have that $W^H(\omega)\in B^H_{\Phi_2, \infty}([0,1];\mathbb{R})$ holds for a.a.\ $\omega\in\Omega$ and therefore, we can apply \autoref{thm:lift} with $\mathbb{W}^H_{s,t}=\tfrac12 (\delta W^H_{s,t})^2$ and $F^H_{s,t}=0$ to obtain
\begin{equation*}
	\mathbf{W}^H(\omega)=(W^H(\omega),\mathbb{W}^H(\omega))\in\mathscr{B}^{H}_{\Phi_2, \infty}\quad\text{for a.a.\ } \omega\in\Omega.
\end{equation*}
Finally, such regularity is also obtained for non-Gaussian processes. For example, if $n\in\mathbb{N}$ and $\beta_1,\beta_2\in \mathbb{R}$ are real numbers such that 
\[
	0<\beta_1 + \frac{n}{2}(\beta_2-1)+1<1, \quad
	1-\frac{1}{n}<\beta_2<1,\quad \mbox{and}\quad 
	\frac{1}{3}<\kappa = \beta_1 + \frac{n}{2}(\beta_2-1)+1\leq \frac{1}{2},
\]
then it holds for the fractionally filtered Hermite process $Z^{\beta_1,\beta_2,n}$ (see \cite[p.318]{CouOnd24} and also \cite[Theorem 3.27]{BaiTaq14} for the definition) that $Z^{\beta_1,\beta_2,n}(\omega) \in B_{\Phi_{2/n},\infty}^\kappa([0,1];\mathbb{R})$ for a.a.\ $\omega\in\Omega$ (see \cite[Corollary 4.2]{CouOnd24}). As such, these processes can be also lifted via \autoref{thm:lift} with $\mathbb{Z}^{\beta_1,\beta_2,n}_{s,t}=\frac{1}{2}(\delta Z^{\beta_1,\beta_2,n}_{s,t})^2$ and $F_{s,t}^{\beta_1,\beta_2,n}=0$ to a rough path that satisfies
\[
	\mathbf{Z}^{\beta_1,\beta_2,n}(\omega) = (Z^{\beta_1,\beta_2,n}(\omega),\mathbb{Z}^{\beta_1,\beta_2,n}(\omega))\in\mathscr{B}^{\kappa}_{\Phi_{2/n}, \infty}\quad\text{for a.a.\ } \omega\in\Omega.
\]
	
\subsection{Controlled rough paths}

In this section, we define a Besov-Orlicz variant of the controlled rough path from \cite{Gub04}, establish some basic estimates, and study the composition of controlled rough paths with sufficiently smooth functions.

For the rest of this section, let $n, m, k\in\mathbb{N}$ be fixed.

\begin{definition}
\label{def:controlled_RP}
	Given $X\in\BOone([0,T]; \R^n)$ and $Y\in\BOone([0,T]; \Rmk)$, we say that $Y$ is \emph{controlled} by $X$ if there is $Y'\in\BOone([0,T]; \Lcal(\R^n; \Rmk))$ such that the remainder $R^{Y}$ given by
	\[
		\delta Y_{s,t} = Y'_s\delta X_{s,t} + R^{Y}_{s, t}, \quad (s,t) \in \triangle^2[0,T],
	\]
	satisfies $R^{Y} \in \BOtwo([0,T]; \Rmk)$. We denote the space of all controlled rough paths, i.e.\ the space of all such pairs $(Y,Y')$, by $\CRPs([0,T]; \Rmk)$. For $(Y,Y')\in\CRPs([0,T];\Rmk)$, we define
	\[
		[(Y,Y')]_{\CRPs([0,T];\Rmk)} = [Y']_{\BOone([0,T];\Lcal(\R^n; \Rmk))} + \norm{R^Y}{\BOtwo([0,T];\Rmk)}.
	\]
	Additionally, for $\tilde{X}\in\BOone([0,T];\R^n)$ and $(\tilde{Y},\tilde{Y}')\in\CRPtld([0,T];\Rmk)$, we denote
	\begin{multline*}
		\dXXon{T}((Y,Y'),(\tilde{Y},\tilde{Y}'))\\
		= d_{\BOone([0,T];\Lcal(\R^n;\Rmk))}(Y',\tilde{Y}') + d_{\BOtwo([0,T];\Rmk)}(R^Y,R^{\tilde{Y}}).
	\end{multline*}
\end{definition}

If there is no risk of ambiguity, we will often omit the domains and codomains and write, e.g., only $\CRPs$ and $\dXX$ instead of $\CRPs([0,T]; \Rmk)$ and $\dXXon{T}$, respectively. We follow a similar convention for other spaces as well. If $X=\tilde{X}$, we use the notation $\dX=\dXX$. Let us note that we follow the convention from \cite{FriHai20} in the sense that the coefficients of $\CRPs$ correspond to the regularity of the remainder $R^Y$. In the often referenced paper \cite{FriSee22}, the coefficients correspond to the regularity of the underlying path $X$ instead.

For brevity, we will often omit the norm-related subscripts and write only, e.g., $|X_t|$, $|Y_t|$, $|Y'_t|$, and $|R^Y_{s,t}|$ instead of $|X_t|_{\R^n}$, $|Y_t|_{\Rmk}$, $|Y'_t|_{\Lcal(\R^n;\Rmk)}$, and $|R^Y_{s,t}|_{\Rmk}$, respectively.

\begin{lemma}\label{lem:deltaRYdeltaRYdiffleqM}
	Let $X\in\BOone([0,T]; \R^n)$, $(Y,Y')\in\CRPs([0,T];\Rmk)$, and $p \in (\frac{1}{\alpha}, \infty)$. Then, for any $(s,u,t)\in\triangle^3([0,T])$, it holds that
	\[
		\lvert\delta R^Y_{s,u,t}\rvert_{\Rmk}\lsmT{\alpha,\beta,p,q}[Y']_{\BOone}[X]_{\BOone}\left(((u-s)\wedge (t-u)^\frac{1}{2}((u-s)\vee(t-u))^\frac{1}{2}\right)^{2(\alpha-\frac{1}{p})}
	\]
	and
	\begin{multline*}
		\lvert\delta(R^Y-R^{\tilde{Y}})_{s,u,t}\rvert_{\Rmk} \lsmT{\alpha,\beta,p,q}\left([Y'-\tilde{Y}']_{\BOone}[X]_{\BOone} + [\tilde{Y}']_{\BOone}[X-\tilde{X}]_{\BOone}\right)
		\\
		\cdot \left(((u-s)\wedge (t-u)^\frac{1}{2}((u-s)\vee(t-u))^\frac{1}{2}\right)^{2(\alpha-\frac{1}{p})}.
	\end{multline*}
\end{lemma}

\begin{proof}
	From the definition of $R^Y$, we deduce
	\begin{align*}
		\delta R^Y_{s,u,t}
		&=-Y'_s\delta X_{u,t} + Y'_u\delta X_{u,t} = (\delta Y'_{s,u})(\delta X_{u,t}), \quad (s,u,t) \in \triangle^3[0, T].
	\end{align*}
	The first claim can be therefore established analogously to \autoref{lem:deltaXdeltaXdiffleqM}. The second claim follows similarly from
	\[
		\delta(R^Y-R^{\tilde{Y}})_{s,u,t} = (\delta(Y'-\tilde{Y}')_{s,u})(\delta X_{u,t}) + (\delta \tilde{Y}'_{s,u})(\delta(X-\tilde{X})_{u,t}), \quad (s,u,t) \in \triangle^3[0, T].
	\]
\end{proof}

\begin{lemma}\label{lem:odhady}
	Let $p \in (\frac1\alpha, \infty)$, $X, \tilde{X}\in\BOone([0,T];\R^n)$, $(Y,Y')\in\CRPs([0,T];\Rmk)$, and $(\tilde{Y},\tilde{Y}')\in\CRPtld([0,T];\Rmk)$. Then
	\begin{align*}
		\norm{R^Y}{C^{2\alpha - 2/p; 2}}&\lsmT{\alpha, \beta, p, q} \norm{R^Y}{\BOtwo} + [Y']_{\BOone}[X]_{\BOone},
		\\
		\norm{R^Y-R^{\tilde{Y}}}{C^{2\alpha - 2/p; 2}}&\lsmT{\alpha, \beta, p, q} \norm{R^Y-R^{\tilde{Y}}}{\BOtwo}
		\\
		&\hphantom{\lsmT{\alpha, \beta, p, q} } + [Y'-\tilde{Y}']_{\BOone}[X]_{\BOone} + [\tilde{Y}']_{\BOone}[X-\tilde{X}]_{\BOone},
	\end{align*}
	and
	\begin{align*}
		[Y]_{\BOone}&\lsmT{\alpha, \beta, p, q} \lvert Y'_0\rvert_{\mathcal{L}(\R^n;\Rmk)}[X]_{\BOone}
		\\
		&\hphantom{\lsmT{\alpha, \beta, p, q} } +(T^{\alpha-\frac{1}{p}}\vee T^{\alpha-\frac{3}{p}})\left([Y']_{\BOone}[X]_{\BOone} + \norm{R^Y}{\BOtwo}\right),
		\\
		[Y-\tilde{Y}]_{\BOone}&\lsmT{\alpha, \beta, p, q} \lvert Y'_0-\tilde{Y}'_0\rvert_{\mathcal{L}(\R^n; \Rmk)}[X]_{\BOone} + \lvert \tilde{Y}'_0\rvert_{\mathcal{L}(\R^n;\Rmk)}[X-\tilde{X}]_{\BOone}
		\\
		&\hphantom{\lsmT{\alpha, \beta, p, q} } + (T^{\alpha-\frac{1}{p}}\vee T^{\alpha-\frac{3}{p}})\left([Y'-\tilde{Y}']_{\BOone}[X]_{\BOone} + [\tilde{Y}']_{\BOone}[X-\tilde{X}]_{\BOone} \right)
		\\
		&\hphantom{\lsmT{\alpha, \beta, p, q} } + (T^{\alpha-\frac{1}{p}}\vee T^{\alpha-\frac{3}{p}}) \norm{R^Y-R^{\tilde{Y}}}{\BOtwo}.
	\end{align*}
\end{lemma}

\begin{proof}
	The first two estimates are direct consequences of \cite[Lemma 5.2]{FriSee22} by the embeddings $\BOone \hook \Bone$ and $\BOd{2} \hook \Besov{\alpha;2}{p}{q}$ from \autoref{thm:BOintoB} and \autoref{thm:emb_BO_B_two_param}, respectively.
	
	To establish the third estimate, let us first derive an auxiliary estimate of $R^Y$. By the interpolation inequality from \autoref{lem:interpolation} (with $\delta=2\alpha-\frac{2}{p}$, $\frac{\beta}{2}$ in place of $\beta$, and $\beta'=\beta$), Young's inequality, the first estimate in this lemma and the inclusion $\BO{2\alpha;2}{\beta/2}{q/2} \subseteq \BO{2\alpha;2}{\beta/2}{q}$ from \autoref{prop:emb_BOq1q2} (with $M=[Y']_{\BOone}[X]_{\BOone}$ from \autoref{lem:deltaRYdeltaRYdiffleqM}), we deduce
	\begin{align*}
		\norm{R^Y}{\BOd{2}} &\lesssim_{\alpha,\beta,p,q} T^{\alpha-\frac{1}{p}}\left(\norm{R^Y}{C^{2\alpha-2/p;2}} + \norm{R^Y}{\BO{2\alpha;2}{\beta/2}{q}}\right)
		\\
		&\lesssim_{\alpha,\beta,p,q}^T T^{\alpha-\frac{1}{p}} \left( \norm{R^Y}{\BOtwo} + [Y']_{\BOone}[X]_{\BOone} + [Y']_{\BOone}[X]_{\BOone} T^{-\frac2p}\right)
		\\
		&\lesssim (T^{\alpha-\frac{1}{p}}\vee T^{\alpha-\frac{3}{p}})\left( \norm{R^Y}{\BOtwo} + [Y']_{\BOone}[X]_{\BOone}\right).
	\end{align*}
	Now, the estimate from \autoref{lem:affinesub} yields
	\begin{align*}
		[Y]_{\BOone}&\lesssim_q \norm{Y'}{L^\infty}[X]_{\BOone} + \norm{R^Y}{\BOd{2}}
		\\
		&\lsmT{\alpha,\beta,p,q}\left(\lvert Y'_0\rvert + T^{\alpha-\frac{1}{p}}[Y']_{\BOone} \right) [X]_{\BOone} + \norm{R^Y}{\BOd{2}},
	\end{align*}
	which leads to the desired bound by the auxiliary estimate above.
	
	The remaining estimate follows from
	\begin{equation*}
		\delta(Y-\tilde{Y})_{s,t} = (Y'-\tilde{Y}')_s \delta X_{s,t} + \tilde{Y}'_s \delta(X-\tilde{X})_{s,t} + R^Y_{s,t} - R^{\tilde{Y}}_{s,t}, \quad (s,t) \in \triangle^2[0,T],
	\end{equation*}
	in a similar manner.
\end{proof}

Next, we establish that a composition of a controlled rough path with a sufficiently regular function is a controlled rough path. A similar result was established for Besov-type spaces in \cite[Proposition 5.2]{FriSee22} with the assumption $f \in C^{2, \nu}_b$ for some $\nu \in (0, 1)$. Below, we provide a proof with the assumption  $f \in C^3_b$. The same assumption was also considered in the Sobolev setting of \cite[Lemma 3.7]{LiuProTei2021}.

\begin{theorem}\label{thm:c2comp}
	Let $X\in\BOone([0,T];\R^n)$ and $(Y,Y')\in\CRPs([0,T];\Rmk)$. For $f\in C_b^2(\Rmk, \Rmn)$, define $f(Y)'_\cdot = ((Df)(Y_\cdot))Y'_\cdot$. Then $(f(Y),f(Y)')\in\CRPs([0,T];\Rmn)$ and the bounds
	\begin{equation*}
		[f(Y)]_{\BOone}\leq\norm{Df}{L^\infty}[Y]_{\BOone}
	\end{equation*}
	and
	\begin{multline*}
		[(f(Y),f(Y)')]_{\CRPs} \lsmT{\alpha, \beta, q}\norm{f}{C_b^2}\left(1+[X]_{\BOone}\right)\\
		\cdot \left(\left(\lvert Y'_0\rvert + [(Y,Y')]_{\CRPs}\right) \vee \left(\lvert Y'_0\rvert + [(Y,Y')]_{\CRPs}\right)^2\right)
	\end{multline*}
	hold. Moreover, let $\tilde{X}\in\BOone([0,T];\R^n)$ and $(\tilde{Y},\tilde{Y}')\in\CRPtld([0,T];\Rmk)$ be such that
	\begin{equation*}
		\left(\lvert Y'_0\rvert + [(Y,Y')]_{\CRPs}\right) \vee \left( \lvert \tilde{Y}'_0\rvert + [(\tilde{Y},\tilde{Y}')]_{\CRPtld} \right) \vee [X]_{\BOone} \vee [\tilde{X}]_{\BOone} \leq M
	\end{equation*}
	for some $M \in(0,\infty)$ and let $f\in C_b^3(\Rmk, \Rmn)$, then
	\begin{multline}
		\label{eq:fYprimeminusfYtildeprime}
		[f(Y)' - f(\tilde{Y})']_{\BOone} \lsmTM{\alpha, \beta, q}\norm{f}{C_b^3}
		\\
		\cdot \left(\lvert Y_0-\tilde{Y}_0\rvert + \lvert Y'_0-\tilde{Y}'_0\rvert + [X-\tilde{X}]_{\BOone} + \dXX((Y,Y'),(\tilde{Y},\tilde{Y}'))\right)
	\end{multline}
	and
	\begin{multline*}
		\norm{R^{f(Y)} - R^{f(\tilde{Y})}}{B^{2\alpha;2}_{\Phi_{\beta / 2},q / 2}} \lsmTM{\alpha, \beta, q} \norm{f}{C_b^{3}}
		\\
		\cdot \left(\lvert Y_0 - \tilde{Y}_0\rvert + \lvert Y'_0 - \tilde{Y}'_0\rvert + [X-\tilde{X}]_{\BOone} + \dXX((Y,Y'),(\tilde{Y},\tilde{Y}'))\right).
	\end{multline*}
\end{theorem}

\begin{proof}
	The proof runs similarly as the proof of \cite[Proposition 5.2]{FriSee22} with minor differences such as the use of \autoref{lem:odhady} in place of \cite[Lemma 5.2]{FriSee22} with a fixed choice of $p \in(\frac3\alpha,\infty)$ where required. We provide a more detailed proof of \eqref{eq:fYprimeminusfYtildeprime}.
	
	Set $\tilde{Z}=f(\tilde{Y})$ and $\tilde{Z}'=Df(\tilde{Y})\tilde{Y}'$, then the (quazi-)triangle inequality implies
	\begin{align*}
		[f(Y)'-f(\tilde{Y}')]_{\BOone}&= [Df(Y)Y' - Df(\tilde{Y})\tilde{Y}']_{\BOone}
		\\
		&\lesssim_q [Df(Y)(Y'-\tilde{Y}')]_{\BOone} + [(Df(Y)-Df(\tilde{Y}))\tilde{Y}']_{\BOone}.
	\end{align*}
	Regarding the first term, since $Df$ is Lipschitz, we observe
	\begin{align*}
		\lvert \delta (Df(Y)(Y'-\tilde{Y}'))_{s,t}\rvert &\leq \lvert Df(Y_t) \delta(Y'-\tilde{Y}')_{s,t}\rvert + \lvert \delta(Df(Y))_{s,t}(Y'_s-\tilde{Y}'_s)\rvert
		\\
		&\leq\norm{Df}{L^\infty}\lvert \delta(Y'-\tilde{Y}')_{s,t}\rvert + \norm{D^2f}{L^\infty} \lvert \delta Y_{s,t} \rvert \lvert Y'_s-\tilde{Y}'_s\rvert
	\end{align*}
	for all $(s,t)\in\triangle^2[0,T]$. From \autoref{lem:affinesub}, it follows
	\begin{align*}
		[Df(Y)(Y'-\tilde{Y}')]_{\BOone} &\lesssim_q \norm{Df}{L^\infty} [Y'-\tilde{Y}']_{\BOone} + \norm{D^2f}{L^\infty} [Y]_{\BOone} \norm{Y'-\tilde{Y}'}{L^\infty}
		\\
		&\lesssim_{q}^{T,M} \norm{f}{C^2_b} \left( \lvert Y_0' - \tilde{Y}_0' \rvert + [Y'-\tilde{Y}']_{\BOone} \right).
	\end{align*}
	Similarly, we rewrite the second term as
	\begin{equation*}
		\delta((Df(Y) - Df(\tilde{Y}))\tilde{Y}')_{s,t} = (Df(Y_t)-Df(\tilde{Y}_t))(\delta \tilde{Y}'_{s,t}) + \delta(Df(Y) - Df(\tilde{Y}))_{s,t}\tilde{Y}'_s
	\end{equation*}
	for all $(s,t)\in\triangle^2[0,T]$. By the Lipschitz continuity of $Df$, we have
	\[
		\lvert \delta((Df(Y) - Df(\tilde{Y}))\tilde{Y}')_{s,t} \rvert \leq \norm{D^2f}{L^\infty} \lvert Y_t - \tilde{Y}_t \rvert  \lvert \delta \tilde{Y}'_{s,t} \rvert + \lvert \delta(Df(Y) - Df(\tilde{Y}))_{s,t} \rvert \norm{\tilde{Y}'}{L^\infty}
	\]	
	for all $(s,t)\in\triangle^2[0,T]$. Then, \autoref{lem:affinesub} and the choice of $M$ yield
	\begin{align*}
		&[(Df(Y) - Df(\tilde{Y}))\tilde{Y}']_{\BOone}
		\\
		&\quad \lesssim_q \norm{D^2f}{L^\infty} \norm{Y-\tilde{Y}}{L^\infty} [\tilde{Y}']_{\BOone} + [Df(Y) - Df(\tilde{Y})]_{\BOone} \norm{\tilde{Y}'}{L^\infty}
		\\
		&\quad \lesssim_q^{T,M} \norm{D^2f}{L^\infty} \left( \lvert Y_0 - \tilde{Y}_0 \rvert + [Y-\tilde{Y}]_{\BOone} \right) + [Df(Y) - Df(\tilde{Y})]_{\BOone}.
	\end{align*}	
	Since it holds that
	\begin{align*}
		\delta(Df(Y) - Df(\tilde{Y}))_{s,t} &= \int_0^1 D^2f(\theta Y_t+(1-\theta)\tilde{Y}_t) \delta(Y-\tilde{Y})_{s,t}\dd\theta
		\\
		&\hphantom{= } \quad - \int_0^1 (D^2f(\theta Y_t+(1-\theta)\tilde{Y}_t)-D^2f(\theta Y_s+(1-\theta)\tilde{Y}_s))(Y_s-\tilde{Y}_s)\dd\theta,
	\end{align*}
	we may use the Lipschitz continuity of $D^2 f$ to deduce
	\begin{align*}
		\lvert \delta(Df(Y) - Df(\tilde{Y}))_{s,t} \rvert \leq \norm{D^2 f}{L^\infty} \lvert \delta(Y-\tilde{Y})_{s,t} \rvert + \norm{D^3 f}{L^\infty} \left( \lvert \delta Y_{s,t} \rvert + \lvert \delta \tilde{Y}_{s,t} \rvert \right) \norm{Y-\tilde{Y}}{L^\infty}
	\end{align*}		
	for all $(s,t)\in\triangle^2[0,T]$. Thus, we obtain, by \autoref{lem:affinesub}, the estimate
	\begin{align*}
		[Df(Y) - Df(\tilde{Y})]_{\BOone} &\lesssim_{q} \norm{D^2 f}{L^\infty} [Y-\tilde{Y}]_{\BOone} + \norm{D^3 f}{L^\infty} \left( [Y]_{\BOone} + [\tilde{Y}]_{\BOone} \right) \norm{Y-\tilde{Y}}{L^\infty}\\
		&\lesssim_{q}^{T,M} \norm{f}{C^3_b} \left( [Y-\tilde{Y}]_{\BOone} + \lvert Y_0 - \tilde{Y}_0 \rvert + [Y-\tilde{Y}]_{\BOone} \right),
	\end{align*}
	which finishes the proof of \eqref{eq:fYprimeminusfYtildeprime} by appealing to \autoref{lem:odhady}.
\end{proof}

\subsection{Sewing lemma}

By a partition of $[0,1]$, we understand a finite sequence $\pi = \{t_i\}_{i=0}^N$ such that $t_0=0$, $t_N=1$ and $t_i<t_{i+1}$ for all $i=0,\dots, N-1$. The norm of the partition $\pi = \{t_i\}_{i=0}^N$ is defined by $|\pi| = \max\{ t_{i+1} - t_i \mid 0 \leq i \leq N-1 \}$. A sequence of partitions $\{\pi^n\}_{n=0}^\infty$ is said to have \emph{vanishing norms} if $|\pi^n| \to 0$ as $n \to \infty$.

For a mapping $\varXi\in L^{0;2}([0,T];\R^m)$ and $\pi$, a partition of $[0, 1]$, we define the partial sum $I^\pi\varXi:\triangle^2[0,T]\to\R^m$ by
\begin{equation*}
	(I^\pi\varXi)_{s,t}=\sum_{i=0}^{N-1} \varXi_{s + t_i(t-s), s+t_{i+1}(t-s)}.
\end{equation*}

\begin{theorem}[Sewing Lemma]\label{thm:sewing}
	Assume that $\varXi\in L^{0;2}([0,T];\R^m)$. Let $\alpha\in(0,1),\gamma\in(1,\infty)$, $\beta_1,\beta_2\in(0,\infty)$, and $q_1,q_2\in(0,\infty]$ be such that
	\[
	   \varXi \in \BOSone([0, T];\R^m) \quad \text{and} \quad \delta\varXi \in \BOSthree([0, T];\R^m).
	\]
	Then there exists $\Iscr=\Iscr(\varXi) \in L^0([0,T];\R^m)$ such that the mapping $\varXi\mapsto\Iscr(\varXi)$ is linear, the convergence
	\begin{equation}\label{eq:sewingConvergence}
		\lim_{N\to\infty}\norm{\delta\Iscr-I^{\pi_N}\varXi}{\BOSr}=0
	\end{equation}
	holds for an arbitrary sequence $\{\pi_N\}_{N=1}^\infty$ of partitions of $[0,1]$ with vanishing norms and any $r$ satisfying $r\geq q_2$ and $r>\frac{1}{\gamma}$, and such that there is the estimate
	\begin{equation}\label{eq:sewingEstimate}
		\norm{\delta\Iscr - \varXi}{\BOStwo}\lesssim_{\gamma,q_2}\norm{\delta\varXi}{\BOSthree}.
	\end{equation}
	Moreover, if there exists $\rho\in(0,\frac{1}{2}]$ and $M\in(0,\infty)$ such that
	\begin{equation}\label{eq:sewingDeltavarXi}
		\lvert\delta\varXi_{s,u,t}\rvert\leq M\left(\left(\left(u-s\right)\wedge\left(t-u\right)\right)^\rho\left(\left(u-s\right)\vee\left(t-u\right)\right)^{1-\rho}\right)^{\gamma-\frac{1}{p}}
	\end{equation}
	holds for any $(s,u,t)\in\triangle^3[0,T]$ and some $p\in(\frac{1}{\gamma},\infty)$, then $\Iscr$ is continuous on $[0,T]$, the bound
	\begin{equation}\label{eq:sewingIfirstEst}
		[\Iscr]_{\BO{\alpha}{\beta_1\wedge\beta_2}{q_1\vee q_2}}\lesssim_{\alpha, \gamma, \beta_1, \beta_2, p, q_1, q_2, T} \norm{\varXi}{\BOSone} + \norm{\delta\varXi}{\BOStwo} + M
	\end{equation}
	holds and $\Iscr$ is uniquely determined by prescribing $\Iscr_0 = \Theta$ for some $\Theta \in \R^m$. If, in addition, $\beta_2<\beta_1$, then, for any $p\in(\frac{1}{\alpha}\vee\frac{1}{\gamma-\alpha},\infty)$, we have
	\begin{equation}\label{aln:sewingLemmaLastEsts1}
		\norm{\delta\Iscr-\varXi}{C^{\gamma-1/p;2}}\lsmT{\gamma, p, q_2, \rho}\norm{\delta\varXi}{\BOSthree} + M
	\end{equation}
	and
	\begin{equation}\label{aln:sewingLemmaLastEsts2}
		[\Iscr]_{\BO{\alpha}{\beta_1}{q_1\vee q_2}}\lsmT{\alpha, \gamma, \beta_1, \beta_2, p, q_1, q_2, \rho} \norm{\varXi}{\BOSone} + T^{\gamma - \frac{1}{p} - \alpha + \frac{\beta_2}{p\beta_1}}\left(\norm{\delta\varXi}{\BOSthree} + M\right) + T^{-\frac{1}{p}}M.
	\end{equation}
\end{theorem}

\begin{proof}
	The existence of $\Iscr \in L^0([0,T];\R^m)$ satisfying convergence \eqref{eq:sewingConvergence} and bound \eqref{eq:sewingEstimate} is proved similarly as in \cite[Theorem 3.1]{FriSee22} since we can replace the Besov(-type) norms by Besov-Orlicz(-type) norms in a straightforward manner. The linearity of the map $\varXi \mapsto I\varXi$ can be established from the linearity of the partial sums $I^{\pi_n}(a\varXi+b\tilde{\varXi}) = aI^{\pi_n}\varXi + bI^{\pi_n}\tilde{\varXi}$ holding for any $a,b\in\R$ and $\tilde{\varXi}$ of the same regularity as $\varXi$ by \eqref{eq:sewingConvergence}.

    Assume that \eqref{eq:sewingDeltavarXi} holds. Bound \eqref{eq:sewingIfirstEst} follows from $[\Iscr]_{\BO{\alpha}{\beta_1\wedge\beta_2}{q_1\vee q_2}} = \norm{\delta \Iscr}{\BO{\alpha;2}{\beta_1\wedge\beta_2}{q_1\vee q_2}}$ by the triangle inequality, \autoref{prop:trivemb}, \eqref{eq:sewingEstimate}, and \autoref{prop:emb_BOq1q2} with \eqref{eq:sewingDeltavarXi}. Since $\delta \Iscr \in \BO{\alpha;2}{\beta_1\wedge\beta_2}{q_1\vee q_2}$ and $\delta(\delta \Iscr)=0$, $\delta \Iscr$ has a continuous and bounded modification by \autoref{thm:emb_BO_B_two_param} and, similarly, so does $\varXi$ by \eqref{eq:sewingDeltavarXi}. The continuity of $\delta \Iscr$ also implies the continuity of $\Iscr$.
	
	Let us now establish that the choice of $\Iscr_0$ identifies the map $\Iscr$ uniquely. Let $\tilde{\Iscr} \in C([0, T], \R^m)$ be such that \eqref{eq:sewingEstimate} holds with $\tilde{\Iscr}$ in place of $\Iscr$ and let $\tilde{\Iscr}_0 = \Theta$. Since
	\[
		|\delta(\Iscr - \tilde{\Iscr})_{r,r+h}| \leq |(\delta\Iscr - \varXi)_{r,r+h}| + |(\delta \tilde{\Iscr}-\varXi)_{r,r+h}|	
	\]
	holds for all $(r,r+h) \in \triangle^2[0,T]$, we may directly estimate
	\[
		[\Iscr - \tilde{\Iscr}]_{B^{\gamma}_{\Phi_{\beta_2},q_2}} = \norm{\delta(\Iscr - \tilde{\Iscr})}{\BOStwo} \lesssim_q \norm{\Iscr - \varXi}{\BOStwo} + \norm{\tilde{\Iscr}-\varXi}{\BOStwo} < \infty.
	\]
	From the embeddings in \autoref{thm:BOintoB} and the standard embeddings of Besov spaces from, e.g., \cite[Sections 3.2.4 and 3.3.1]{Tri83}, we obtain
	\[
		B^{\gamma}_{\Phi_{\beta_2},q_2} \hook B^{\gamma}_{p,q_2} \hook B^{\gamma}_{p,\infty} \hook B^{\gamma-\frac1p}_{\infty,\infty} = C^{\gamma-\frac1p}
	\]
	for any $p \in [1, \infty)$. Hence, for $p$ sufficiently large so that $\gamma-\frac1p>1$, it follows that $\Iscr-\tilde{\Iscr}$ is constant on $[0, T]$. By $\Iscr_0 - \tilde{\Iscr}_0 = 0$, it follows that $\Iscr = \tilde{\Iscr}$.
	
	Bounds \eqref{aln:sewingLemmaLastEsts1} and \eqref{aln:sewingLemmaLastEsts2} are proved similarly to \cite[Theorem 3.3 part (a)]{FriSee22}.
	Choose any $p\in(\frac{1}{\alpha}\vee\frac{1}{\gamma-\alpha},\infty)$. \autoref{thm:emb_BO_B_two_param} and bound \eqref{eq:sewingEstimate} implies
	\[
	   \norm{\delta\Iscr-\varXi}{C^{\gamma-1/p}} \lsmT{\gamma, \beta_2, p, \rho}\norm{\delta\Iscr-\varXi}{\BO{\gamma;2}{\beta_2}{q_2}} + M \lesssim_{\gamma, q_2}\norm{\delta\varXi}{\BO{\gamma;3}{\beta_2}{q_2}} + M.
	\]
	Finally, \autoref{lem:interpolation} (with $\delta=\gamma-\frac{1}{p}, \beta'=\beta_1, \beta=\beta_2$, and $q=q_2$) and the above inequality yield
	\begin{align*}
		\norm{\delta\Iscr-\varXi}{\BO{\alpha;2}{\beta_1}{q_2}}&\lesssim_{\alpha,\gamma,\beta_1,\beta_2,p,q_2} T^{\gamma-\alpha-\frac{1}{p}+\frac{\beta_2}{p\beta_1}}\normto{\delta\Iscr-\varXi}{C^{\gamma-1/p;2}}{1-\frac{\beta_2}{\beta_1}}\normto{\delta\Iscr-\varXi}{\BO{\gamma;2}{\beta_2}{q_2}}{\frac{\beta_2}{\beta_1}}
		\\
		&\lsmT{\gamma,\beta_2,p,q_2,\rho}T^{\gamma-\alpha-\frac{1}{p}+\frac{\beta_2}{p\beta_1}}\left(\norm{\delta\varXi}{\BO{\gamma;3}{\beta_2}{q_2}} + M\right)^{1-\frac{\beta_2}{\beta_1}}\left(\norm{\delta\varXi}{\BO{\gamma;3}{\beta_2}{q_2}} + M\right)^{\frac{\beta_2}{\beta_1}}
		\\
		&=T^{\gamma-\alpha-\frac{1}{p}+\frac{\beta_2}{p\beta_1}}\left(\norm{\delta\varXi}{\BO{\gamma;3}{\beta_2}{q_2}} + M\right).
	\end{align*}
	Hence, as we have $\delta(\delta\Iscr-\varXi) = -\delta\varXi$, \autoref{thm:emb_BO_B_two_param} gives
	\begin{align*}
		[\Iscr]_{\BO{\alpha}{\beta_1}{q_1\vee q_2}} &\leq \norm{\varXi}{\BO{\alpha;2}{\beta_1}{q_1\vee q_2}} + \norm{\delta\Iscr-\varXi}{\BO{\alpha;2}{\beta_1}{q_1\vee q_2}}
		\\
		&\lsmT{\alpha,\gamma, q_1, q_2} \norm{\varXi}{\BO{\alpha;2}{\beta_1}{q_1}} + T^{-\frac{1}{p}}M+ \norm{\delta\Iscr-\varXi}{\BO{\alpha;2}{\beta_1}{q_2}} + T^{-\frac{1}{p}}M
		\\
		&\lsmT{\alpha,\gamma,\beta_1,\beta_2,p,q_2,\rho}\norm{\varXi}{\BO{\alpha;2}{\beta_1}{q_1}} + T^{\gamma-\alpha-\frac{1}{p}+\frac{\beta_2}{p\beta_1}}\left(\norm{\delta\varXi}{\BO{\gamma;3}{\beta_2}{q_2}} + M\right) + T^{-\frac{1}{p}}M.
	\end{align*}
\end{proof}

\subsection{Rough integral}

We are now ready to construct the rough integral using the sewing lemma from the previous section. Recall that $\alpha\in(\frac{1}{3}, \frac{1}{2}]$, $\beta\in(0,\infty)$, $q\in(0, \infty]$, let $k=n$, and let $\X\in\BOrp([0,T];\R^n)$, $(Y,Y')\in\CRPs([0,T];\Rmn)$ be fixed.
We define the approximation of the rough integral $\varXi = \varXi_{\X}(Y, Y'): \triangle^2[0, T] \to \R^m$ by
\begin{equation}
	\label{eq:Xi_rough_definition}
	\varXi_{s,t} = Y_s (X_t - X_s) + Y'_s\mathbb{X}_{s, t}, \quad (s,t) \in \triangle^2[0, T].
\end{equation}

Notice that $\varXi\in L^{0;d}([0,T];\R^{m})$. Indeed, we can suppose that $\varXi$ is continuous because it has a continuous modification due to \autoref{thm:BOintoB}, \autoref{thm:emb_BO_B_two_param}, and \autoref{lem:deltaXdeltaXdiffleqM}.
We begin with a simple auxiliary bound.

\begin{lemma}\label{lem:ustuodhady}
	For $(s,u,t) \in \triangle^3[0, T]$ and $\theta\in(0, \infty)$, it holds
	\begin{equation*}
		\lvert u-s\rvert^{2\theta} \lvert t-u\rvert^{\theta} + \lvert u-s\rvert^{\theta} \lvert t-u\rvert^{2\theta} \leq 2 \left(((u-s)\wedge(t-u))^\frac{1}{3}((u-s)\vee(t-u))^\frac{2}{3}\right)^{3\theta}.
	\end{equation*}
\end{lemma}

\begin{proof}
	If $u-s \leq t-u$, then
	\begin{align*}
		\lvert u-s\rvert^{2\theta} \lvert t-u\rvert^{\theta} + \lvert u-s\rvert^{\theta} \lvert t-u\rvert^{2\theta} &= \lvert t-u\rvert^{2\theta}\lvert u-s \rvert^{\theta} \left(\frac{\lvert u-s \rvert^{\theta}}{\lvert t-u\rvert^{\theta}} + 1 \right)
		\\
		&\leq 2 \lvert t-u\rvert^{2\theta}\lvert u-s\rvert^{\theta}
		\\
		&=2 \left(((u-s)\wedge(t-u))^\frac{1}{3} ((u-s)\vee(t-u))^\frac{2}{3}\right)^{3\theta}.
	\end{align*}
	The case $u-s>t-u$ follows similarly.
\end{proof}

Next, we summarize the basic properties of the rough approximation $\varXi$.

\begin{lemma}\label{lem:rozpisDeltaVarXi}
	Let $\varXi: \triangle^2[0, T] \to \R^m$ be as in \eqref{eq:Xi_rough_definition}. Then, for $(s,u,t) \in \triangle^3[0, T]$, it holds
	\[
		\delta\varXi_{s,u,t} = -R^Y_{s,u}\delta X_{u,t} - \delta Y'_{s,u}\Xbb_{u,t}
	\]
	and, for $p \in (\frac1\alpha, \infty)$, we have
	\begin{multline*}
		\lvert\delta\varXi_{s,u,t}\rvert_{\R^m}\lsmT{\alpha,\beta,p,q}\left(\norm{R^Y}{\BOtwo}[X]_{\BOone} + [Y']_{\BOone}\vertiii{\X}{\BOrp}{2}\right)
		\\
		\cdot \left(((u-s)\wedge(t-u))^\frac{1}{3}((u-s)\vee(t-u))^\frac{2}{3}\right)^{3(\alpha-\frac1p)}.
	\end{multline*}
\end{lemma}

\begin{proof}
	The identity follows directly from Chen's relations and the definition of $R^Y$. The bound can be then established in a straightforward manner from the $k(\alpha-1/p)$-H\"older continuity of $R^Y$, $X$, $Y'$, and $\Xbb$ for appropriate $k = 1, 2$ with the embeddings into the Besov-Orlicz and the Besov-Orlicz-type spaces from \autoref{lem:odhady}, \autoref{thm:BOintoB}, and \autoref{lem:preHolderRPEmbs} together with \autoref{lem:ustuodhady}.
\end{proof}

\begin{lemma}\label{lem:rozpisDeltaVarXidiff}
	Let $\mathbf{\tilde{X}} = (\tilde{X}, \tilde{\Xbb}) \in \BOrp([0,T];\R^n)$ and $(\tilde{Y},\tilde{Y}')\in\CRPtld([0,T];\Rmn)$ and let $\tilde{\varXi} = \tilde{\varXi}_{\mathbf{\tilde{X
	}}}(\tilde{Y}, \tilde{Y}')$ be defined by \eqref{eq:Xi_rough_definition}. Assuming that
	\begin{multline*}
		\left(\lvert Y_0\rvert_{\Rmn} + \lvert Y'_0\rvert_{\mathcal{L}(\R^n;\Rmn)} + [(Y,Y')]_{\CRPs}\right) \vee \vertiii{\X}{\BOrp}{}
		\\
		\vee \left(\lvert \tilde{Y}_0\rvert_{\Rmn} + \lvert \tilde{Y}'_0\rvert_{\mathcal{L}(\R^n;\Rmn)} + [(\tilde{Y},\tilde{Y}')]_{\CRPtld} \right) \vee \vertiii{\mathbf{\tilde{X}}}{\BOrp}{} \leq M
	\end{multline*}
	holds for some $M\in[0,\infty)$, then, for all $(s,u,t) \in \triangle^3[0, T]$, we have the estimate
	\begin{multline*}
		\lvert\delta(\varXi-\tilde{\varXi})_{s,u,t}\rvert \lsmTM{\alpha,\beta,p,q} \left(((u-s)\wedge(t-u))^\frac{1}{3}((u-s)\vee(t-u))^\frac{2}{3}\right)^{3(\alpha-\frac{1}{p})}
		\\
		\cdot \left(\norm{R^Y-R^{\tilde{Y}}}{\BOtwo} + [Y'-\tilde{Y}']_{\BOone} + [X-\tilde{X}]_{\BOone} + \norm{\Xbb-\tilde{\Xbb}}{\BOtwo}\right).
	\end{multline*}
\end{lemma}

\begin{proof}
	Since it is readily seen that 
	\[
		\delta(\varXi-\tilde{\varXi})_{s,u,t} = (R^{\tilde{Y}}-R^Y)_{s,u}\delta X_{u,t} + R^{\tilde{Y}}_{s,u}\delta(\tilde{X} - X)_{u,t}+ \delta (\tilde{Y}'-Y')_{s,u}\Xbb_{u,t} + \delta\tilde{Y}'_{s,u}(\tilde{\Xbb}-\Xbb)_{u,t}
	\]
	holds for all $(s,u,t)\in\triangle^3[0,T]$, we may use the H\"older continuity similarly as in the proof of \autoref{lem:rozpisDeltaVarXi} to deduce
	\begin{multline*}
		\lvert\delta(\varXi-\tilde{\varXi})_{s,u,t}\rvert_{\R^m}	\lsmM{} \left(\lvert u-s\rvert^{2\alpha-\frac{2}{p}}\lvert t-u\rvert^{\alpha-\frac{1}{p}} + \lvert u-s\rvert^{\alpha-\frac{1}{p}}\lvert t-u\rvert^{2\alpha-\frac{2}{p}}\right)
		\\
		\cdot \left(\norm{R^Y-R^{\tilde{Y}}}{C^{2\alpha-2/p;2}} + [X-\tilde{X}]_{C^{\alpha-1/p}} + [Y'-\tilde{Y}']_{C^{\alpha-1/p}} + \norm{\Xbb-\tilde{\Xbb}}{C^{2\alpha-2/p;2}}\right).
	\end{multline*}
	The rest of the claim follows from the estimates in \autoref{lem:odhady}, \autoref{thm:BOintoB} and \autoref{thm:emb_BO_B_two_param} and the bound in \autoref{lem:ustuodhady}.
\end{proof}

We are now ready to establish the existence of the rough integral that corresponds to the local approximation $\varXi$.

\begin{theorem}\label{thm:gubinelli}
    The mapping
    \begin{equation*}
	   \I = \I_{\X}(Y,Y') : [0, T] \to \R^m, \quad t \mapsto\Iscr(\varXi_{\X}(Y,Y'))_t,
    \end{equation*}
	where $\Iscr(\varXi_{\X}(Y,Y'))$ is the map from \autoref{thm:sewing} and $\varXi_{\X}(Y,Y')$ is defined in~\eqref{eq:Xi_rough_definition}, is well-defined. Moreover, $\I_0=0$, the mapping
	\[
		\mathscr{I} : \CRPs([0,T];\Rmn) \to \CRPs([0,T];\R^m), \quad (\bar{Y},\bar{Y}') \mapsto \left(\I_{\X}(\bar{Y},\bar{Y}'),\bar{Y}\right)
	\]	
	is linear, the limit
	\begin{equation}
		\label{eq:gubinelli_zero_limit}
		\lim\limits_{N\to\infty}\norm{\delta\I - I^{\pi_N}\varXi}{\BOthreer} = 0
	\end{equation}
	holds for any $\bar{q}>\frac{1}{3\alpha}$, $\bar{q}\geq\frac{q}{3}$, and any sequence $\{\pi_N\}_{N=1}^\infty$ of partitions of $[0, 1]$ with vanishing norm,
	and the bounds
	\begin{equation}
		\label{eq:gubinelli_dI-Xi_bound}
		\norm{\delta\I - \varXi}{\BOthree} \lsmT{\alpha, \beta, q} \norm{\delta\varXi}{\BO{3\alpha;3}{\beta/3}{q/3}} \lsmT{q} \norm{R^Y}{\BOtwo}[X]_{\BOone} + [Y']_{\BOone}\vertiii{\X}{\BOrp}{}
	\end{equation}
	and
	\begin{multline}
		\label{eq:gubinelli_I_besov_seminorm}
		[\I]_{\BOone} \lesssim_{\alpha,\beta,q,T} \left( \lvert Y_0\rvert_{\Rmn} + [Y]_{\BOone} + \norm{R^Y}{\BOtwo} \right) [X]_{\BOone}\\
		+ \left(\lvert Y'_0\rvert_{\Lcal(\R^n;\Rmn)} + [Y']_{\BOone} \right)\vertiii{\X}{\BOrp}{2}
	\end{multline}
	hold. Additionally, for any $p \in (\frac{3}{\alpha}, \infty)$, it holds that
	\begin{multline}
		\label{eq:gubinelli_RI}
		\norm{R^{\I}}{\BOtwo}\lsmT{\alpha, \beta, p, q}\lvert Y'_0\rvert_{\Lcal(\R^n;\Rmn)}\norm{\Xbb}{\BOtwo} 
		\\
		+ T^{\alpha-\frac{1}{p}} \left(\vertiii{\X}{\BOrp}{}\vee\vertiii{\X}{\BOrp}{2}\right)[(Y,Y')]_{\CRPs}.
	\end{multline}
\end{theorem}

\begin{proof}
	With, e.g., $p=4/\alpha$, it holds that
	\begin{equation}\label{eq:gubvarxi}
		\norm{\varXi}{\BOd{2}}\lsmT{\alpha, \beta, q} \norm{Y}{L^\infty} [X]_{\BOone} + \norm{Y'}{L^\infty} \vertiii{\X}{\BOrp}{2}.
	\end{equation}
	Recalling the definition of $\varXi$ in \eqref{eq:Xi_rough_definition}, a straightforward application of the H\"older-type inequality in Orlicz spaces from \autoref{thm:orliczholder} (with $\beta_1=3$, $\beta_2=\frac{3}{2}$, and $\beta/3$ in place of $\beta$) and \autoref{lem:luxp} leads to
	\begin{equation}\label{eq:gubinelideltavarXi}
		\norm{\delta\varXi}{\BO{3\alpha;3}{\beta/3}{q/3}} \lsmT{q} \norm{R^Y}{\BOtwo}[X]_{\BOone} + [Y']_{\BOone}\norm{\Xbb}{\BOtwo}.
	\end{equation}
	Since \eqref{eq:gubvarxi}, \eqref{eq:gubinelideltavarXi}, and the bound in \autoref{lem:rozpisDeltaVarXi} verify the assumptions of the sewing lemma in \autoref{thm:sewing}, there exists $\I\in\OspF([0,T];\R^m)$ that satisfies $\delta\I - \varXi \in\BOthree([0,T];\R^m)$ and \eqref{eq:gubinelli_dI-Xi_bound} and, for which \eqref{eq:gubinelli_zero_limit} holds for any sequence $\{\pi_N\}_{N=1}^\infty$ of partitions of $[0,1]$ with vanishing norms and for all $\bar{q}>\frac{1}{3\alpha}$ and $\bar{q}\geq\frac{q}{3}$. Estimate \eqref{eq:gubinelli_I_besov_seminorm} follows from \eqref{aln:sewingLemmaLastEsts2} in \autoref{thm:sewing} by \eqref{eq:gubvarxi}, \eqref{eq:gubinelideltavarXi}, and by the bound in \autoref{lem:rozpisDeltaVarXi} with \autoref{lem:affinesub} and the choice, e.g., $p=2/\alpha$.
	
	It remains to establish \eqref{eq:gubinelli_RI}. Let $p \in (\frac3\alpha, \infty)$. By setting $\I' = Y$, we observe that
	\[
		\delta\I - \varXi = \delta\I - \I'\delta X - Y'\Xbb = R^\I - Y'\Xbb
	\]
	holds and thus,
	\[
		\norm{R^\I}{\BOtwo}\leq \norm{Y'\Xbb}{\BOtwo} + \norm{\delta\I - \varXi}{\BOtwo}
	\]
	is obtained. Regarding the first term in the above inequality, by \autoref{lem:affinesub} and the definition of norms on $\CRPs$ and $\BOrp$, we deduce
		\begin{align*}
		\norm{Y'\Xbb}{\BOtwo}&\leq\norm{Y'}{L^\infty}\norm{\Xbb}{\BOtwo}
		\\
		&\lsmT{\alpha,\beta, p, q} \lvert Y'_0\rvert_{\Lcal(\R^n;\Rmn)}\norm{\Xbb}{\BOtwo} + T^{\alpha-\frac{1}{p}}[(Y,Y')]_{\CRPs}\vertiii{\X}{\BOrp}{2}.
	\end{align*}
	For the second term, let us denote $\mathcal{R} \varXi = \delta \I - \varXi$ for brevity. Since $\delta \mathcal{R} \varXi = \delta(\delta \I - \varXi) = - \delta \varXi$, the inclusion $B^{3\alpha; 2}_{\Phi_{\beta/3}, q/3} \subseteq C^{3(\alpha-1/p)}$ from \autoref{thm:emb_BO_B_two_param} (with the bound from \autoref{lem:rozpisDeltaVarXi}) and \eqref{eq:gubinelideltavarXi} yield
	\begin{equation}\label{eq:gubineliRinC}
	\begin{split}
		\norm{\mathcal{R} \varXi}{C^{3(\alpha-\frac1p)}}&\lsmT{\alpha,p,q}\norm{\delta\varXi}{\BO{3\alpha;3}{\beta/3}{q/3}} + \norm{R^Y}{\BOtwo}[X]_{\BOone} + [Y']_{\BOone}\vertiii{\X}{\BOrp}{2}
		\\
		&\leq\left(\vertiii{\X}{\BOrp}{}\vee\vertiii{\X}{\BOrp}{2}\right)[(Y,Y')]_{\CRPs}.
	\end{split}
	\end{equation}
	By appealing to the interpolation inequality from \autoref{lem:interpolation} (with $3\delta = \alpha-\frac{3}{p}$, $\beta' = \frac{\beta}{2}$, and $\frac{\beta}{3}$ in place of $\beta$),
	the Young inequality, the inclusion $\BO{3\alpha;2}{\beta/3}{q/3} \subseteq \BO{3\alpha;2}{\beta/3}{q/2}$ from \autoref{prop:emb_BOq1q2}, and to estimate \eqref{eq:gubinelli_dI-Xi_bound}, we obtain
	\begin{align*}
		\norm{\mathcal{R} \varXi}{\BOtwo}&\lsmT{\alpha,\beta,p,q} T^{2\alpha-\frac{1}{p}}\left(\norm{\mathcal{R} \varXi}{C^{3\alpha-1/p;2}} + \norm{\mathcal{R} \varXi}{\BO{3\alpha;2}{\beta/3}{q/2}} \right)
		\\
		&\lsmT{\alpha,q}T^{2\alpha-\frac{1}{p}}\left(\norm{\mathcal{R} \varXi}{C^{3\alpha-1/p;2}} + \norm{\mathcal{R} \varXi}{\BOthree} \right)
		\\
		&\hphantom{\lsmT{\alpha,q} } \ + T^{2\alpha-\frac{4}{p}}\left(\norm{R^Y}{\BOtwo}[X]_{\BOone} + [Y']_{\BOone}\norm{\Xbb}{\BOtwo}\right)
		\\
		&\lesssim_{\alpha,q}T^{2\alpha-\frac{1}{p}}\left(\norm{\mathcal{R} \varXi}{C^{3\alpha-1/p;2}} + \norm{\delta\varXi}{\BOthree}\right)
		\\
		&\hphantom{\lsmT{\alpha,q} } \ + T^{2\alpha-\frac{4}{p}}\left(\norm{R^Y}{\BOtwo}[X]_{\BOone} + [Y']_{\BOone}\norm{\Xbb}{\BOtwo}\right).
	\end{align*}
	Hence, bounds  \eqref{eq:gubineliRinC} and \eqref{eq:gubinelideltavarXi} imply
	\[
		\norm{\mathcal{R} \varXi}{\BOtwo}\lsmT{\alpha,\beta,p,q} T^{2\alpha-\frac{4}{p}} \left(\vertiii{\X}{\BOrp}{}\vee\vertiii{\X}{\BOrp}{2}\right)[(Y,Y')]_{\CRPs}
	\]
	and the proof of \eqref{eq:gubinelli_RI} is concluded by gathering the estimates above.
\end{proof}

\begin{definition}\label{def:roughintegral}
	For $\X \in \BOrp([0,T];\R^n)$ and $(Y,Y')\in\CRPs([0,T];\Rmn)$, we define the \emph{rough integral of $Y$ with respect to the rough path $\X$} by
	\[
		\int_0^\cdot Y_s \dd\X_s: [0, T] \to \R^m, \quad \int_0^\cdot Y_s \dd\X_s = \I_{\X}(Y, Y')_\cdot,
	\]
	where $\I_{\X}(Y,Y')$ is defined in \autoref{thm:gubinelli}.
\end{definition}

We conclude this section with a continuity estimate for the rough integral.

\begin{theorem}\label{thm:estimatefordX}
	Let $\X, \mathbf{\tilde{X}}\in\BOrp([0,T];\R^n)$, $(Y,Y')\in\CRPs([0,T];\Rmn)$, and $(\tilde{Y},\tilde{Y}')\in\CRPtld([0,T];\Rmn)$. Assume that
	\begin{multline*}
		\left(\lvert Y_0\rvert_{\Rmn} + \lvert Y'_0\rvert_{\mathcal{L}(\R^n;\Rmn)} + [(Y,Y')]_{\CRPs}\right) \vee \vertiii{\X}{\BOrp}{}
		\\
		\vee \left(\lvert \tilde{Y}_0\rvert_{\Rmn} + \lvert \tilde{Y}'_0\rvert_{\mathcal{L}(\R^n;\Rmn)} + [(\tilde{Y},\tilde{Y}')]_{\CRPtld} \right) \vee \vertiii{\mathbf{\tilde{X}}}{\BOrp}{} \leq M
	\end{multline*}
	holds for some $M \in(0,\infty)$. Let $\tilde{\varXi}=\varXi_{\mathbf{\tilde{X}}}(\tilde{Y},\tilde{Y}')$ and $\tilde{\I}=\I_{\mathbf{\tilde{X}}}(\tilde{Y},\tilde{Y}')$, where $\varXi$ and $\I$ are defined in \eqref{eq:Xi_rough_definition} and \autoref{thm:gubinelli}, respectively. Then, for any $p \in ( \frac{4}{\alpha}, \infty)$, it holds that 
	\begin{multline*}
		\dXX(\left( \I, Y\right), (\tilde{\I}, \tilde{Y})) \lsmTM{\alpha,\beta,p,q} \lvert Y_0-\tilde{Y}_0\rvert^{\frac{q}{2}\wedge 1}_{\Rmn} + \lvert Y'_0-\tilde{Y}'_0\rvert_{\Lcal(\R^n;\Rmn)}^{\frac{q}{2}\wedge 1}
		\\
		+\rho_{\BOrp}(\X,\mathbf{\tilde{X}})^{\frac{q}{2}\wedge 1} + (T^{\alpha-\frac{4}{p}})^{\frac{q}{2}\wedge 1} \left( \normto{R^Y - R^{\tilde{Y}}}{\BOtwo}{\frac{q}{2}\wedge 1} + [Y' - \tilde{Y}']^{\frac{q}{2}\wedge 1}_{\BOone} \right).
	\end{multline*}
\end{theorem}

\begin{proof}
	The proof is similar to the proof of \cite[Theorem 5.5]{FriSee22}. However, since the resulting bounds are different and since these bounds are important for our main result, we include the details.

	We will often implicitly use the inequality $x \leq (1 \vee M) x^{q \wedge 1}$ that holds for all $q \in (0, \infty)$ and $x \in [0, M]$ and the following two estimates.
	From \autoref{lem:HolderRPEmbs}, we observe 
	\[
	\norm{\Xbb}{\BOd{2}} \lsmT{\alpha,\beta,p,q} (T^{\alpha-\frac{1}{p}}\vee T^{\alpha-\frac{2}{p}}) \vertiii{\X}{\BOrp}{2} \lsmTM{} T^{\alpha-\frac{2}{p}}
	\]
	and
	\begin{align*}
		\norm{\Xbb-\tilde{\Xbb}}{\BOd{2}} &\lsmT{\alpha,\beta,p,q} \sum_{i=1}^2\left(\vertiii{\X}{\BOrp}{}\vee\vertiii{\mathbf{\tilde{X}}}{\BOrp}{}\right)^{2-i}\norm{\X^{(i)}-\mathbf{\tilde{X}}^{(i)}}{\BO{i\alpha;2}{\beta/i}{q/i}}
		\\
		&\lsmTM{q} \rho_{\BOrp}(\X,\mathbf{\tilde{X}}).
	\end{align*}
	Also note that by \autoref{lem:affinesub}, \autoref{lem:odhady}, and the choice of $M$, it holds
	\begin{align*}
		\norm{Y}{L^\infty}&\lsmT{\alpha,\beta,p,q}\lvert Y_0\rvert + T^{\alpha-\frac{1}{p}}[Y]_{\BOone}
		\\
		&\lsmTM{\alpha,\beta,p,q} \lvert Y_0\rvert + T^{\alpha-\frac1p} \left( \lvert Y'_0\rvert + (T^{\alpha-\frac{1}{p}}\vee T^{\alpha-\frac{3}{p}}) \left([Y']_{\BOone} + \norm{R^Y}{\BOtwo}\right) \right)
		\\
		&\lsmTM{} 1.
	\end{align*}
	Similarly, \autoref{lem:affinesub} implies $\norm{Y'}{L^\infty} \lsmTM{\alpha,\beta,p,q} 1$. Bounds of the same type obviously hold also for $\tilde{Y}$ and $\tilde{Y}'$. In the rest of the proof, we will use these estimates implicitly as well.
	
	By \autoref{lem:affinesub}, \autoref{lem:odhady}, and by the inequality $\alpha-\frac3p \leq \alpha-\frac1p$, we estimate
	\begin{align*}
		&\dXX((\I,Y),(\tilde{\I},\tilde{Y}'))
		\\
		&\quad \lesssim_q \norm{Y-\tilde{Y}}{\Osp} + [Y-\tilde{Y}]_{\BOone}^{q\wedge 1} + \normto{R^\I-R^{\tilde{\I}}}{\BOtwo}{\frac{q}{2}\wedge 1}
		\\
		&\quad \lsmTM{\alpha,\beta,p,q} \lvert Y_0-\tilde{Y}_0\rvert + [Y-\tilde{Y}]_{\BOone}^{q\wedge 1} + \normto{R^\I-R^{\tilde{\I}}}{\BOtwo}{\frac{q}{2}\wedge 1}
		\\
		&\quad \lsmTM{\alpha,\beta,p,q} \lvert Y_0-\tilde{Y}_0\rvert + \lvert Y_0'-\tilde{Y}_0'\rvert^{q \wedge 1} + \rho_{\BOrp}(\X,\mathbf{\tilde{X}})^{q \wedge 1}
		\\
		&\quad \hphantom{\lsmTM{\alpha,\beta,p,q} \ } +(T^{\alpha-\frac{3}{p}})^{q\wedge 1} [Y'-\tilde{Y}']_{\BOone}^{q\wedge 1} + \normto{R^\I-R^{\tilde{\I}}}{\BOtwo}{\frac{q}{2}\wedge 1}.
	\end{align*}
	Recalling the equalities $\mathcal{R}\varXi=\delta \I - \varXi = \delta \I - Y\delta X - Y'\Xbb$ and $\mathcal{R}\tilde{\varXi} = \delta\tilde{\I} -\tilde{Y}\delta\tilde{X} - \tilde{Y}'\tilde{\Xbb}$, we can estimate the last term by
	\begin{equation}
		\label{eq:RZdiff}
		\norm{R^\I-R^{\tilde{\I}}}{\BOtwo} \lesssim_q \norm{\mathcal{R}\varXi - \mathcal{R}\tilde{\varXi}}{\BOtwo} + \normto{Y'\Xbb - \tilde{Y}'\tilde{\Xbb}}{\BOtwo}{}.
 	\end{equation}
	
	Before we proceed to the estimates of the first term on the right-hand side of \eqref{eq:RZdiff}, we establish several auxiliary bounds. Similarly as in the proof of \autoref{thm:gubinelli}, we use the bound from \autoref{lem:rozpisDeltaVarXidiff} with \autoref{thm:orliczholder} and \autoref{lem:luxp} to deduce
	\begin{equation}
		\label{eq:rough_int_stab_deltaXi}
		\norm{\delta(\varXi - \tilde{\varXi})}{\BO{3\alpha;3}{\beta/3}{q/3}} \lsmTM{q} \rho_{\BOrp}(\X,\mathbf{\tilde{X}}) + [Y'-\tilde{Y}']_{\BOone} + \norm{R^Y-R^{\tilde{Y}}}{\BOtwo}.
	\end{equation}
	By \autoref{lem:affinesub} and \autoref{lem:odhady} and by observing that $\alpha-\frac1p \leq 2\alpha-\frac4p$ holds by the choice of $p$, we have
	\begin{align*}
		\norm{Y-\tilde{Y}}{L^\infty}&\lsmT{\alpha,\beta,p,q}\lvert Y_0-\tilde{Y}_0\rvert + T^{\alpha-\frac{1}{p}}[Y-\tilde{Y}]_{\BOone}
		\\
		&\lsmTM{\alpha,\beta,p,q} \lvert Y_0-\tilde{Y}_0\rvert + \lvert Y'_0-\tilde{Y}'_0\rvert + \rho_{\BOrp}(\X,\mathbf{\tilde{X}})
		\\
		&\hphantom{\lsmTM{\alpha,\beta,p,q} \ } + T^{\alpha-\frac1p} \left([Y'-\tilde{Y}']_{\BOone} + \norm{R^Y-R^{\tilde{Y}}}{\BOtwo}\right)
	\end{align*}
	and, by \autoref{lem:affinesub}, we also have
	\[
		\norm{Y'-\tilde{Y}'}{L^\infty}\lsmT{\alpha,\beta,p,q}\lvert Y_0'-\tilde{Y}_0'\rvert_{\Lcal(\R^n;\Rmn)} + T^{\alpha-\frac{1}{p}}[Y'-\tilde{Y}']_{\BOone}.			
	\]
	Noticing that $\delta \mathcal{R}(\varXi - \tilde{\varXi}) = - \delta(\varXi - \tilde{\varXi})$, the embedding into H\"older-type functions from  \autoref{thm:emb_BO_B_two_param} (with the bound from \autoref{lem:rozpisDeltaVarXidiff}) and \eqref{eq:rough_int_stab_deltaXi} yield
	\begin{align*}
		\norm{\mathcal{R}(\varXi-\tilde{\varXi})}{C^{3(\alpha-\frac{1}{p});2}} &\lsmTM{\alpha,p,q} \norm{\delta(\varXi-\tilde{\varXi})}{\BO{3\alpha;3}{\beta/3}{q/3}} + \rho_{\BOrp}(\X,\mathbf{\tilde{X}})
		\\
		&\hphantom{\lsmTM{\alpha,p,q} \ } + [Y'-\tilde{Y}']_{\BOone} + \norm{R^Y-R^{\tilde{Y}}}{\BOtwo}
		\\
		&\lsmTM{q} \rho_{\BOrp}(\X,\mathbf{\tilde{X}}) + [Y'-\tilde{Y}']_{\BOone} + \norm{R^Y-R^{\tilde{Y}}}{\BOtwo}.
	\end{align*}
	Finally, we use the interpolation inequality from \autoref{lem:interpolation} (with $\alpha = 2\alpha$, $\gamma= 3\alpha$, $\frac{\beta}{3}$ in place of $\beta$, $\beta'=\frac{\beta}{2}$, and $\frac{q}{2}$ in place of $q$), the Young inequality, the inclusion $ \BOthree \subseteq B^{3\alpha; 2}_{\Phi_{\beta/3}, q/2}$ from \autoref{prop:emb_BOq1q2}, estimate \eqref{eq:sewingEstimate} from \autoref{thm:sewing} (applied to $\varXi - \tilde{\varXi}$ in place of $\varXi$ thanks to linearity) and the bound above with estimate \eqref{eq:rough_int_stab_deltaXi} to get
	\begin{align*}
		&\norm{\mathcal{R}(\varXi - \tilde{\varXi})}{\BOtwo} \lesssim_{\alpha,\beta,p,q} T^{\alpha-\frac{1}{p}}\left(\norm{\mathcal{R}(\varXi - \tilde{\varXi})}{C^{3(\alpha-\frac{1}{p});2}} + \norm{\mathcal{R}(\varXi - \tilde{\varXi})}{\BO{3\alpha;3}{\beta/3}{q/2}}\right)
		\\
		&\quad \lsmT{\alpha, p} T^{\alpha-\frac{1}{p}} \left( \norm{\mathcal{R}(\varXi - \tilde{\varXi})}{C^{3(\alpha-\frac{1}{p});2}} + \norm{\mathcal{R}(\varXi - \tilde{\varXi})}{\BOthree} \right)
		\\
		&\quad \hphantom{\lsmT{\alpha, p} \ } + T^{\alpha-\frac1p} T^{-\frac3p} \left( [Y'-\tilde{Y}']_{\BOone} + \norm{R^Y-R^{\tilde{Y}}}{\BOtwo} \right) + \rho_{\BOrp}(\X,\mathbf{\tilde{X}})
		\\
		&\quad \lsmTM{\alpha,\beta,p,q} \rho_{\BOrp}(\X,\mathbf{\tilde{X}}) + T^{\alpha-\frac4p} \left( [Y'-\tilde{Y}']_{\BOone} + \norm{R^Y-R^{\tilde{Y}}}{\BOtwo} \right).
	\end{align*}
	
	It remains to estimate the second term on the right-hand side of \eqref{eq:RZdiff}. By the (quazi-)triangle inequality, the choice of $M$ and the bounds above, we deduce
	\begin{align*}
		\norm{Y'\Xbb - \tilde{Y}'\tilde{\Xbb}}{\BOtwo} &\lesssim_q \norm{Y'}{L^\infty}\norm{\Xbb-\tilde{\Xbb}}{\BOtwo} + \norm{Y'-\tilde{Y}'}{L^\infty}\norm{\tilde{\Xbb}}{\BOtwo}
		\\
		&\lsmTM{\alpha,\beta,p,q} \rho_{\BOrp}(\X,\mathbf{\tilde{X}}) + \lvert Y_0'-\tilde{Y}_0'\rvert + T^{\alpha-\frac{4}{p}}[Y'-\tilde{Y}']_{\BOone}.
	\end{align*}
	The proof is concluded by collecting the bounds above and straightforward estimates using the inequality $x \leq (1 \vee M) x^{q \wedge 1}$ from the beginning of the proof.
\end{proof}

\section{Rough differential equations with Besov-Orlicz signals}
\label{sec:RDEs}

In this section, we state and prove our main result -- the existence and uniqueness of solutions to rough differential equations driven by exponential Besov-Orlicz signals. We fix $\alpha\in(\frac{1}{3},\frac{1}{2}]$, $\beta\in(0,\infty)$, and $q\in(0,\infty]$ for the whole section. Let also $\X\in\BOrp([0,T];\R^n)$ and $f\in C_b^3(\R^m,\Rmn)$. 

First, we define the notion of a solution to a rough differential equation.

\begin{definition}
	We say that $(Y,Y')\in\CRPs([0,T];\R^m)$ is a \emph{solution} to the rough differential equation (RDE)
	\begin{equation}\label{eq:rdeobs}
		\dd Y_t = f(Y_t) \dd \X_t,\quad Y_0=y\in\R^m,
	\end{equation}
	if
	\[Y_t = y + \int_0^t f(Y_s)\dd\X_s,\]
	where the integral $\int_0^t f(Y_s)\dd\X_s$ is the rough integral from \autoref{def:roughintegral}, holds for all $t\in[0,T]$.
\end{definition}

We note that the rough integral $\int_0^\cdot f(Y_s)\dd\X_s$ in the above definition is well-defined. Indeed, \autoref{thm:c2comp} guarantees that
\begin{equation}\label{eq:rdeRIreg1}
	(f(Y),f(Y)')\in\CRPs([0,T];\Rmn).
\end{equation}
Therefore, \autoref{thm:gubinelli} ensures that
\begin{equation}\label{eq:rdeRIreg2}
	y+\int_0^\cdot f(Y)\dd\X\in\BOone([0,T];\R^m)
\end{equation}
and that
\begin{equation}\label{eq:rdeRIreg3}
	R^{y+\int_0^\cdot f(Y)\dd\X}\in\BOtwo([0,T];\R^m).
\end{equation}

We are now ready to state the main result.

\begin{theorem}\label{thm:existanduniq}
	For any $y\in\R^m$, there exists a unique solution $(Y,Y')\in\CRPs([0,T];\R^m)$ to RDE \eqref{eq:rdeobs}.
\end{theorem}

Before we proceed to the proof, we include an auxiliary completeness result used in the proof of \autoref{thm:existanduniq}.

\begin{lemma}\label{lem:ycompl}
	The space
	\[\YscrT=\{(Y,Y')\in\CRPs([0,T];\R^m)\,|\,Y_0=y, Y'_0=f(y)\}\]
	equipped with the metric $\dX$ from \autoref{def:controlled_RP} is complete.
\end{lemma}

\begin{proof}
	Let $\{(\Yn,\Ynp)\}_{n\in\mathbb{N}}$ be a Cauchy sequence in $\YscrT$. By \autoref{def:controlled_RP}, we need to find a subsequence converging to a point in $\YscrT$ in the metrics on the spaces $\BOone$ and $\BOtwo$ defined in \eqref{eq:1stMetric} and \eqref{eq:2ndMetric}, respectively. By the completeness of the affine subspaces from \autoref{lem:affinesub} and the bounds in \autoref{lem:odhady}, we may find $Y\in\BOone([0,T]; \R^m)$, $Y'\in\BOone([0,T]; \mathcal{L}(\R^n;\R^m))$, and $R^Y\in\BOtwo([0,T];\R^m)$ such that $Y_0=y$, $Y'_0=f(y)$, and
	\begin{equation*}
		\lim\limits_{n\to\infty}[\Yn - Y]_{\BOone}=0,\quad
		\lim\limits_{n\to\infty}[\Ynp - Y']_{\BOone}=0,\quad
		\lim\limits_{n\to\infty}\norm{R^{\Yn}-R}{\BOtwo}=0.
	\end{equation*}
	Moreover, the estimate from \autoref{lem:affinesub} yields
	\begin{equation}\label{eq:orliczconv}
		\norm{\Ynp-Y'}{L^{\infty}}\vee\norm{\Ynp-Y'}{\Osp}\lesssim_{\alpha,\beta,q,T}[\Ynp-Y']_{\BOone}.
	\end{equation}
	Choose $\varepsilon\in(0,\alpha)$ arbitrary. Then, using $R^{\Yn}=\delta\Yn-\Ynp\delta X$ and the embedding $B^{2\alpha;2}_{\Phi_{\beta/2},q} \hook B^{\alpha-\varepsilon;2}_{\Phi_{\beta/2},q}$ from \autoref{prop:trivemb}, we obtain
	\begin{align*}
		&\norm{R-(\delta Y-Y'\delta X)}{\BOeps}
		\\
		&\quad\lesssim_q\norm{R-R^{\Yn}}{\BOeps} + \norm{(\delta\Yn-\Ynp\delta X)-(\delta Y-Y'\delta X)}{\BOeps}
		\\
		&\quad\lsmT{\alpha,\beta,q,\varepsilon}\norm{R-R^{\Yn}}{\BOtwo} + \norm{\delta(\Yn-Y)}{\BO{\alpha-\varepsilon;2}{\beta}{q/2}} + \norm{\Ynp-Y'}{L^{\infty}}\norm{\delta X}{\BO{\alpha-\varepsilon;2}{\beta}{q/2}}
		\\
		&\quad=\norm{R-R^{\Yn}}{\BOtwo} + [\Yn-Y]_{\BO{\alpha-\varepsilon}{\beta}{q/2}} + \norm{\Ynp-Y'}{L^{\infty}}\left [X\right]_{\BO{\alpha-\varepsilon}{\beta}{q/2}}.
	\end{align*}
	Hence, the embedding $B^{\alpha}_{\Phi_\beta,q} \hook B^{\alpha-\varepsilon}_{\Phi_\beta,q/2}$ from \autoref{cor:boeps} and estimate \eqref{eq:orliczconv} yield
	\begin{multline*}
		\norm{R-(\delta Y-Y'\delta X)}{\BOeps}
		\\
		\lesssim_{\alpha,\beta,q,\varepsilon,T}
		\norm{R-R^{\Yn}}{\BOtwo} + [\Yn-Y]_{\BOone} + [\Ynp-Y']_{\BOone}\left [X\right]_{\BOone}.
	\end{multline*}
	As $n\to\infty$, the right-hand side approaches zero. Consequently, $(Y,Y')\in\CRPs([0,T];\R^m)$, implying that $(Y,Y')\in\YscrT$. Finally, estimate~\eqref{eq:orliczconv} yields
	\[\lim\limits_{n\to\infty}\dX((\Yn,\Ynp),(Y,Y'))=0,\]
	which concludes the proof.
\end{proof}

\begin{proof}[Proof of \autoref{thm:existanduniq}]
	We aim to use the Banach fixed-point argument to prove \autoref{thm:existanduniq}. To this end, let
	\begin{equation*}
		\bar{Y}_t=y+f(y)X_t\quad\text{and}\quad \bar{Y}'_t \equiv f(y)\quad\text{for}\ t\in[0,T],
	\end{equation*}
	and note that $R^{\bar{Y}}\equiv0$ on $\triangle^2[0,T]$. For $M \geq 0$, we denote the closed ball in $\YscrT$ of radius $M$ centred in $(\bar{Y},\bar{Y}')$ by
	\begin{equation*}
		\BYM = \{(Y,Y')\in\YscrT\,|\,\dX((Y,Y'),(\bar{Y},\bar{Y}'))\leq M\}.
	\end{equation*}
	By \autoref{lem:ycompl}, $\BYM$ is complete. Let $\Zscr:\YscrT\to\YscrT$ be defined by
	\begin{equation*}
		\Zscr(Y,Y')=\left(y+\int_0^\cdot f(Y)\dd\X,f(Y)\right).
	\end{equation*}
	The map $\Zscr$ is well-defined thanks to regularities \eqref{eq:rdeRIreg1}, \eqref{eq:rdeRIreg2}, and \eqref{eq:rdeRIreg3}, and the initial conditions $(\int_0^\cdot f(Y)\dd\X)_0=0$ and $f(Y_0)=f(y)$. We will find finite positive constants $M_0$ large enough and $T_0$ small enough, both independent of $y$, such that
	\begin{enumerate}[label=\roman*)]
		\item $\Zscr(\BYMMT{M_0}{\YscrTzero})\subseteq\BYMMT{M_0}{\YscrTzero}$, i.e. \label{p:onball}
		\[
			\norm{f(Y)-\bar{Y}'}{\Osp} + [f(Y)-\bar{Y}']^{q\wedge 1}_{\BOone} + \normto{R^{y+\int_0^\cdot f(Y)\dd\X}-R^{\bar{Y}}}{\BOtwo}{\frac{q}{2}\wedge 1}\leq M_0
		\]
		for all $(Y,Y')\in\BYMMT{M_0}{\YscrTzero}$; and 
		\item $\Zscr$ is a contraction on $\BYMMT{M_0}{\YscrTzero}$, i.e. \label{p:contr}
		\[\dX(\Zscr(Y,Y'),\Zscr(\tilde{Y},\tilde{Y}'))\leq \tilde{C}\dX((Y,Y'),(\tilde{Y},\tilde{Y}'))\]
		for some $\tilde{C}<1$ and all $(Y,Y'),(\tilde{Y},\tilde{Y}')\in\BYMMT{M_0}{\YscrTzero}$.
	\end{enumerate}
	The Banach fixed-point theorem then yields the existence and uniqueness of a local solution, i.e.\ a solution on $[0, T_0]$, that can be extended to the full interval $[0,T]$ in a finite number of steps.
	
	To establish claim \ref{p:onball} above, fix $M\in(0,\infty)$, $(Y,Y')\in\BYM$, and $p\in(\frac{4}{\alpha},\infty)$. Since
	\[
		[(Y,Y')]_{\CRPs}=[Y']_{\BOone}+\norm{R^Y}{\BOtwo}=[Y'-f(y)]_{\BOone}+\norm{R^Y-R^{\bar{Y}}}{\BOtwo},
	\]
	the estimate
	\begin{align*}
		[(Y,Y')]_{\CRPs} &\leq\dX((Y,Y'),(\bar{Y},\bar{Y}'))^{1 / (q\wedge 1)} + \dX((Y,Y'),(\bar{Y},\bar{Y}'))^{1 / (\frac{q}{2}\wedge 1)}
		\\
		&\leq M^{1 / (q\wedge 1)} + M^{1 / (\frac{q}{2}\wedge 1)} =: F_q(M). \label{eq:RDEthmpoint(i)firstEq}\numberthis
	\end{align*}
	holds. Then, from \autoref{thm:c2comp}, we obtain
	\begin{align*}
		&[(f(Y),f(Y)')]_{\CRPs}
		\\
		&\quad\lsmT{\alpha,\beta,q} \left(\norm{f}{C_b^3}\vee\normto{f}{C_b^3}{2}\right)\left(1\vee\vertiii{\X}{\BOrp}{}\right)\left(1+\left([(Y,Y')]_{\CRPs}\vee[(Y,Y')]_{\CRPs}^2\right)\right)
		\\
		\numberthis
		&\quad\leq \hat{C}_{(T)} \left(1+(F_q(M)\vee F_q(M)^2)\right)),\label{eq:fyfyprimeleqM}
	\end{align*}
	where $\hat{C}_{(T)}=C(\norm{f}{C^3_b([0,T])}, \vertiii{\X}{\BOrp([0,T])}{})$ for some nondecreasing positive function $C$ captures the dependence on the respective norms of $f$ and $\X$ as a function of $T$. Thus, as a function of $T$, $\hat{C}_{(T)}$ is also nondecreasing. Even though the precise value of $\hat{C}_{(T)}$, i.e.\ the particular function $C$, may increase from line to line, it will still retain the form described above. Since the Besov-Orlicz seminorm is invariant to shifts, from the Lipschitz continuity of $f$, it follows that
	\[
		[f(Y)-\bar{Y}']_{\BOone}=[f(Y)-f(y)]_{\BOone}=[f(Y)]_{\BOone}\leq\norm{Df}{L^\infty}[Y]_{\BOone}.
	\]
	By using \autoref{lem:odhady} and estimate \eqref{eq:RDEthmpoint(i)firstEq}, we obtain
	\begin{equation*}
		\begin{split}
			&[f(Y)-\bar{Y}']_{\BOone}
			\\
			&\quad\lsmT{\alpha,\beta,p,q}\norm{Df}{L^\infty}\left(\absh{Y'_0}{\R^{m\times n}}[X]_{\BOone} + T^{\alpha-\frac{3}{p}}\left([Y']_{\BOone}[X]_{\BOone} + \norm{R^Y}{\BOtwo}\right)\right)
			\\
			&\quad\leq\norm{f}{C_b^3}\left(\norm{f}{L^\infty}\vertiii{\X}{\BOrp}{} + T^{\alpha-\frac{3}{p}}\left([Y']_{\BOone}\vertiii{\X}{\BOrp}{} + \norm{R^Y}{\BOtwo}\right)\right)
			\\
			&\quad\leq \hat{C}_{(T)}(1+T^{\alpha-\frac{3}{p}}F_q(M)).
		\end{split}
	\end{equation*}
	Thus, by \autoref{lem:affinesub} and $\bar{Y}'_0=f(y)=f(Y_0)$, we deduce
	\begin{equation}\label{eq:point1_1}
		\norm{f(Y)-\bar{Y}'}{\Osp} + [f(Y)-\bar{Y}']^{q\wedge 1}_{\BOone}\lesssim_{\alpha,\beta,p,q} \hat{C}_{(T)} (1+T^{\alpha-\frac{3}{p}}F_q(M)).
	\end{equation}
	We remark that we could drop the exponent $q\wedge 1$ as $(1+T^{\alpha-\frac{3}{p}}(F_q(M))) \geq 1$. We will use this simplification from now on without explicit notice.
	In a similar manner, since
	\[
		\norm{R^{y+\int_0^\cdot f(Y)\dd\X}-R^{\bar{Y}}}{\BOtwo}=\norm{R^{y+\int_0^\cdot f(Y)\dd\X}}{\BOtwo}=\norm{R^{\int_0^\cdot f(Y)\dd\X}}{\BOtwo},
	\]
	estimate \eqref{eq:gubinelli_RI} from \autoref{thm:gubinelli} for the remainder term $R^{\int_0^\cdot f(Y)\dd\X}$ and estimate \eqref{eq:fyfyprimeleqM} yield
	\begin{equation}\label{eq:point1_2}
		\begin{split}
			&\norm{R^{y+\int_0^\cdot f(Y)\dd\X}-R^{\bar{Y}}}{\BOtwo}
			\\
			&\quad \lsmT{\alpha,\beta,p,q} \absh{f(Y)'_0}{\Lcal(\R^n;\Rmn)}\norm{\Xbb}{\BOtwo}
			\\
			&\quad\hphantom{\lsmT{\alpha,\beta,p,q}} \ + T^{\alpha-\frac{1}{p}}\left(\vertiii{\X}{\BOrp}{}\vee\vertiii{\X}{\BOrp}{2}\right)[(f(Y),f(Y)')]_{\CRPs}
			\\
			&\quad\lsmT{\alpha,\beta,q}\hat{C}_{(T)}(1+T^{\alpha-\frac{1}{p}}(1+(F_q(M)\vee F_q(M)^2))).
		\end{split}
	\end{equation}
	Combining estimates \eqref{eq:point1_1} and \eqref{eq:point1_2}, we observe
	\begin{multline}
		\label{eq:mballest}
		\norm{f(Y)-\bar{Y}'}{\Osp\onT} + [f(Y)-\bar{Y}']^{q\wedge 1}_{\BOone\onT} + \normto{R^{y+\int_0^\cdot f(Y)\dd\X}-R^{\bar{Y}}}{\BOtwo\onT}{\frac{q}{2}\wedge 1}
		\\
		\leq C_1\left(\alpha,\beta,p,q,T\right)\hat{C}_{(T)} \left(1 +T^{\alpha-\frac{3}{p}}\left(1+\left(F_q(M)\vee F_q(M)^2\right)\right)\right)
	\end{multline}
	where $C_1$ is a fixed constant with nondecreasing dependence on $T$. Recalling \eqref{eq:fyfyprimeleqM}, we may assume that $C_1$ and $\hat{C}_{(T)}$ also satisfy
	\begin{equation}\label{eq:candfchoice}
			[(f(Y),f(Y)')]_{\CRPs}\leq C_1\left(\alpha,\beta,p,q,T\right)\hat{C}_{(T)} \left(1+\left(F_q(M)\vee F_q(M)^2\right)\right)
	\end{equation}
	without any loss of generality.
	We set
	\[
		M_0 = \left( 2C_1(\alpha,\beta,p,q,T) \hat{C}_{(T)} \right)
		\vee\vertiii{\X}{\BOrp\onT}{}\vee\norm{f}{C_b^{3}\onT}
	\]
	and we choose $T_1\in(0,T]$ so that
	\begin{equation*}
		T_1^{\alpha-\frac{3}{p}}\leq\frac{1}{1+(F_q(M_0)\vee F_q(M_0)^2)}.
	\end{equation*}
	Hence, by rewriting estimate \eqref{eq:mballest} for $(Y,Y')\in\BYMMT{M_0}{\YscrTone}$ and by recalling the nondecreasing dependence of $\hat{C}_{(T)}$ and $C_1$ on $T$, we obtain
	\begin{align*}
		&\norm{f(Y)-\bar{Y}'}{\Osp([0,T_1])} + [f(Y)-\bar{Y}']_{\BOone([0,T_1])}^{q\wedge 1} + \normto{R^{y+\int_0^\cdot f(Y)\dd\X}-R^{\bar{Y}}}{\BOtwo([0,T_1])}{\frac{q}{2}\wedge 1}
		\\
		&\quad\leq C_1(\alpha,\beta,p,q,T_1)\hat{C}_{(T_1)} \left(1 +T_1^{\alpha-\frac{3}{p}}\left(1+\left(F_q(M_0)\vee F_q(M_0)^2\right)\right)\right)
		\\
		&\quad\leq 2 C_1(\alpha,\beta,p,q,T)\hat{C}_{(T)} \leq M_0,
	\end{align*}
	which completes the proof of claim \ref{p:onball}. Let us note that the claim also holds for any $\hat{T}_1\in(0,T_1]$.
	
	It remains to establish claim \ref{p:contr}. We first observe that, for all $(Y,Y')\in\BYMMT{M_0}{\YscrTone}$, estimate \eqref{eq:fyfyprimeleqM} together with the choice of $C_1$ and~$\hat{C}_{(T)}$ so that \eqref{eq:candfchoice} holds as well, implies
	\begin{multline*}
		\absh{f(Y)_0}{\R^{m\times n}}\vee\absh{f(Y)'_0}{\Lcal(\R^n;\Rmn)}\vee[(f(Y),f(Y)')]_{\CRPs}
		\\
		\lesssim^{T_1}_{\alpha,\beta,q} (F_q(M_0)\vee F_q(M_0)^2)(1\vee F_q(M_0))\left(1+\left(F_q(M_0)\vee F_q(M_0)^2\right)\right) \lesssim^{T_1,M_0}_{\alpha,\beta,p,q} 1.
	\end{multline*}
	Let us fix $(Y,Y'),(\tilde{Y},\tilde{Y}')\in\BYMMT{M_0}{\YscrTone}$. With the bound above, \autoref{thm:estimatefordX} yields
	\begin{multline*}
		\dXon{T_1}(\Zscr(Y,Y'),\Zscr(\tilde{Y},\tilde{Y}'))
		\\
		\lsmToneMzero{\alpha,\beta,p,q}T_1^{(\alpha-\frac{4}{p})(\frac{q}{2}\wedge1)}\left(\normto{R^{f(Y)}-R^{f(\tilde{Y})}}{\BOgamma([0,T_1])}{\frac{q}{2}\wedge1}+[f(Y)'-f(\tilde{Y})']^{\frac{q}{2}\wedge 1}_{\BOone([0,T_1])}\right).
	\end{multline*}
	Hence, from \autoref{thm:c2comp}, we obtain
	\begin{multline}
		\label{eq:contraction.preliminary}
		\dXon{T_1}(\Zscr(Y,Y'),\Zscr(\tilde{Y},\tilde{Y}'))
		\\
		\leq C_2(\alpha,\beta,p,q,M_0,T_1) \bar{C}_{(T_1)} T_1^{(\alpha-\frac{4}{p})(\frac{q}{2}\wedge1)}\dXon{T_1}((Y,Y'),(\tilde{Y},\tilde{Y}')),
	\end{multline}
	where $\bar{C}_{(T_1)}=C(\norm{f}{C^3_b([0,T_1])}, \vertiii{\X}{\BOrp([0,T_1])}{})$ denotes a constant of the same form as $\hat{C}_{(T)}$ and $C_2(\alpha,\beta,p,q,M_0,T_1)$ is a constant nondecreasing in the last argument. Let $T_0\in(0,T_1]$ be such that
	\[
		T_0^{(\alpha-\frac{4}{p})(\frac{q}{2}\wedge1)}<\frac{1}{2C_2(\alpha,\beta,p,q,M_0,T_1)\bar{C}_{(T_1)}}.
	\]
	We may now consider estimate \eqref{eq:contraction.preliminary} on $[0, T_0]$ instead of $[0,T_1]$. By the nondecreasing dependence of both $C_2$ and $\bar{C}_{(T_1)}$ on $T_1$, we deduce
	\begin{align*}
		&\dXon{T_0}(\Zscr(Y,Y'),\Zscr(\tilde{Y},\tilde{Y}'))
		\\
		&\quad\leq C_2(\alpha,\beta,p,q,M_0,T_0) \bar{C}_{(T_0)} T_0^{(\alpha-\frac{4}{p})(\frac{q}{2}\wedge1)}\dXon{T_0}((Y,Y'),(\tilde{Y},\tilde{Y}'))
		\\
		&\quad\leq C_2(\alpha,\beta,p,q,M_0,T_1)\bar{C}_{(T_1)}T_0^{(\alpha-\frac{4}{p})(\frac{q}{2}\wedge1)}\dXon{T_0}((Y,Y'),(\tilde{Y},\tilde{Y}'))
		\\
		&\quad<\frac{1}{2}\dXon{T_0}((Y,Y'),(\tilde{Y},\tilde{Y}')),
	\end{align*}
	which concludes the proof of \autoref{thm:existanduniq}.
\end{proof}

\textbf{Acknowledgements:} 
FH has received funding from the Charles University Grant Agency (GA~UK), project no.\ 178823.
P\v{C} and JS were supported by the Czech Science Foundation project no.\ 22-12790S

\medskip

\textbf{Declarations of interest:} None.


\begin{thebibliography}{10}

\bibitem{BaiTaq14}
S.~Bai and M.~S. Taqqu.
\newblock Generalized {H}ermite processes, discrete chaos and limit theorems.
\newblock {\em Stoch. Proc. Appl.}, 124(4):1710--1739, 2014.

\bibitem{Cie91}
Z.~Ciesielski.
\newblock Modulus of smoothness of the {B}rownian paths in the {$L^p$} norm.
\newblock In {\em Constructive theory of functions}, pages 71--75, Varna,
  Bulgaria, 1991.

\bibitem{Cie93}
Z.~Ciesielski.
\newblock {O}rlicz spaces, spline systems, and {B}rownian motion.
\newblock {\em Constr. Approx.}, 9:191--208, 1993.

\bibitem{CieKerRoy93}
Z.~Ciesielski, G.~Kerkyacharian, and B.~Roynette.
\newblock Quelques espaces fonctionnels associ\'es \`a des processus gaussiens.
\newblock {\em Stud. Math.}, 107(2):171--204, 1993.

\bibitem{CouOnd24}
P.~{\v{C}}oupek and M.~Ondrej\'at.
\newblock {B}esov-{O}rlicz path regularity of non-{G}aussian processes.
\newblock {\em Potential Anal.}, 60:307--339, 2024.

\bibitem{DecUst99}
L.~Decreusefond and A.~S. \"{U}st\"{u}nel.
\newblock Stochastic analysis of the fractional {B}rownian motion.
\newblock {\em Potential Anal.}, 10(2):177--214, 1999.

\bibitem{FriHai20}
P.~K. Friz and M.~Hairer.
\newblock {\em A course on rough paths}.
\newblock Universitext. Springer Cham, 2 edition, 2020.

\bibitem{FriPro18}
P.~K. Friz and D.~J. Pr\"omel.
\newblock Rough path metrics on a {B}esov-{N}ikolskii-type scale.
\newblock {\em Tran. Amer. Math. Soc.}, 379(12):8521--8550, 2018.

\bibitem{FriSee22}
P.~K. Friz and B.~Seeger.
\newblock {B}esov rough path analysis.
\newblock {\em J. Diff. Equ.}, 339:152--231, 2022.

\bibitem{Gub04}
M.~Gubinelli.
\newblock Controlling rough paths.
\newblock {\em J. Funct. Anal.}, 216(1):86--140, 2004.

\bibitem{Lev37}
P.~L\'evy.
\newblock {\em Th\'eorie de l'addition des variables al\'eatoires}, volume xvii
  of {\em Monographies des Probabilit\'es; calcul des probabilit\'es et ses
  applications, publi\'ees sous la direction de E. Borel}.
\newblock Gauthier-Villars, Paris, 1937.

\bibitem{LiuProTei21}
Ch. Liu, D.~J. Pr\"omel, and J.~Teichmann.
\newblock On {S}obolev rough paths.
\newblock {\em J. Math. Anal. Appl.}, 497(1):124876, 2021.

\bibitem{LiuProTei2021}
Ch. Liu, D.~J. Pr\"{o}mel, and J.~Teichmann.
\newblock On {S}obolev rough paths.
\newblock {\em J. Math. Anal. Appl.}, 497(1):Paper No. 124876, 21, 2021.

\bibitem{LiuProTei23}
Ch. Liu, D.~J. Pr\"omel, and J.~Teichmann.
\newblock Optimal extension to {S}obolev rough paths.
\newblock {\em Potential Anal.}, 59:1399--1424, 2023.

\bibitem{LiuProTei23a}
Ch. Liu, D.~J. Pr\"omel, and J.~Teichmann.
\newblock A {S}obolev rough path extension theorem via regularity structures.
\newblock {\em ESAIM: Prob. Stat.}, 27:136--155, 2023.

\bibitem{OndVer20}
M.~Ondrej\'at and M.~C. Veraar.
\newblock On temporal regularity of stochastic convolutions in 2-smooth
  {B}anach spaces.
\newblock {\em Ann. Inst. H. Poincar\'e Probab. Statist.}, 56(3):1792--1808,
  2020.

\bibitem{OndSimKup18}
M.~Ondrej\'at, P.~\v{S}imon, and M.~Kupsa.
\newblock Support of solutions of stochastic differential equations in
  exponential {B}esov--{O}rlicz spaces.
\newblock {\em Stoch. Anal. Appl.}, 36:1037--1052, 2018.

\bibitem{PicKufJohFuc13}
L.~Pick, A.~Kufner, O.~John, and S.~Fu\v{c}\'{\i}k.
\newblock {\em Function spaces. {V}ol. 1}, volume~14 of {\em De Gruyter Series
  in Nonlinear Analysis and Applications}.
\newblock Walter de Gruyter \& Co., Berlin, extended edition, 2013.

\bibitem{ProTra16}
D.~J. Pr\"omel and M.~Trabs.
\newblock Rough differential equations driven by signals in {B}esov spaces.
\newblock {\em J. Diff. Equ.}, 260(6):5202--5249, 2016.

\bibitem{RaoRen91}
M.~M. Rao and Z.~D. Ren.
\newblock {\em Theory of {O}rlicz spaces}.
\newblock M. Dekker, 1991.

\bibitem{Tri83}
H.~Triebel.
\newblock {\em Theory of Function Spaces}.
\newblock Birkh\"auser Basel, 1983.

\bibitem{Tud08}
C.~A. Tudor.
\newblock Analysis of the {R}osenblatt process.
\newblock {\em ESAIM: Prob. Stat.}, 12:230--257, 2008.

\bibitem{Ver09}
M.~C. Veraar.
\newblock Correlation inequalities and applications to vector-valued {G}aussian
  random variables and fractional {B}rownian motion.
\newblock {\em Potential Anal.}, 30(4):341--370, 2009.

\bibitem{Wich21}
J.~Wichmann.
\newblock On temporal regularity for strong solutions to stochastic
  $p$-{L}aplace systems.
\newblock 2021.
\newblock Preprint, \href{https://arxiv.org/abs/2111.09601}{arXiv:2111.09601}.

\end{thebibliography}
\end{document}